\documentclass[a4paper, 11pt]{imsart}
\usepackage[OT1]{fontenc}
\usepackage{amsfonts,amsthm,amssymb,amsmath}
\usepackage{mathrsfs}
\usepackage{xcolor}

\startlocaldefs
\theoremstyle{plain}

\endlocaldefs

\usepackage{psfrag, epsfig, graphicx}
\usepackage{amsfonts,amsmath,latexsym,amssymb}
\usepackage[active]{srcltx}
\usepackage{fullpage}

\newtheorem{lemma}{Lemma}
\newtheorem{theorem}{Theorem}
\newtheorem{proposition}{Proposition}

\newtheorem*{assmN1}{Assumption (N1)}
\newtheorem*{assmN2}{Assumption (N2)}
\newtheorem*{assmN3}{Assumption (N3)}

\newtheorem*{ass4.1}{Assumption $4^\prime$}
\newtheorem*{ass4.2}{Assumption $4^{\prime\prime}$}

\newtheorem{definition}{Definition}

\renewcommand{\kappa}{\varkappa}

\newcommand{\e}{\varepsilon}
\newcommand{\rf}{\mathrm{f}}

\newcommand{\rd}{\mathrm{d}}
\newcommand{\ra}{\mathrm{a}}
\newcommand{\rb}{\mathrm{b}}

\newcommand{\cC}{{\cal C}}

\newcommand{\cH}{{\cal H}}
\newcommand{\cI}{{\cal I}}

\newcommand{\cO}{{\cal O}}
\newcommand{\cP}{{\cal P}}

\newcommand{\cR}{{\cal R}}

\newcommand{\cT}{{\cal T}}
\newcommand{\cU}{{\cal U}}

\newcommand{\cX}{{\cal X}}

\newcommand{\bB}{\mathbb B}

\newcommand{\bE}{\mathbb E}
\newcommand{\bF}{\mathbb F}

\newcommand{\bL}{{\mathbb L}}

\newcommand{\bN}{{\mathbb N}}

\newcommand{\bP}{{\mathbb P}}

\newcommand{\bR}{{\mathbb R}}
\newcommand{\bS}{{\mathbb S}}

\newcommand{\bZ}{{\mathbb Z}}
\newcommand{\mA}{\mathfrak{A}}

\newcommand{\mF}{\mathfrak{F}}
\newcommand{\mH}{\mathfrak{H}}
\newcommand{\mS}{\mathfrak{S}}

\newcommand{\mP}{\mathfrak{P}}

\newcommand{\md}{\mathfrak{d}}

\newcommand{\mT}{\mathfrak{T}}

\newcommand{\mE}{\mathfrak{E}}

\newcommand{\epr}{\hfill\hbox{\hskip 4pt
                \vrule width 5pt height 6pt depth 1.5pt}\vspace{0.5cm}\par}

\begin{document}

\begin{frontmatter}
\title{Structural adaptive deconvolution under $L_p$-losses}
\runtitle{Lp adaptive deconvolution}

\begin{aug}
\author{\fnms{Gilles} \snm{Rebelles}
\ead[label=e1]{rebelles.gilles@neuf.fr}}
\runauthor{G. Rebelles}
\affiliation{ Aix--Marseille Universit\'e }
\address{Institut de Math\'ematique de Marseille \\
Aix-Marseille Universit\'e \\
39 rue F. Joliot-Curie 13453 Marseille, France\\
\printead{e1}\\ }
\end{aug}

\begin{abstract}
In this paper, we address the problem of estimating a multidimensional density $f$ by using indirect observations from the statistical model $Y=X+\varepsilon$. Here, $\varepsilon$ is a measurement error independent of the random vector $X$ of interest, and having a known density with respect to the Lebesgue measure. Our aim is to obtain optimal accuracy of estimation under  $\bL_p$-losses when the error $\varepsilon$ has a characteristic function with a polynomial decay. To achieve this goal, we first construct a kernel estimator of $f$ which is fully data driven. Then, we derive for it an oracle inequality under very mild assumptions on the characteristic function of the error $\varepsilon$. As a consequence, we get minimax adaptive upper bounds over a large scale of anisotropic Nikolskii classes and we prove that our estimator is asymptotically rate optimal when $p\in[2,+\infty]$. Furthermore, our estimation procedure adapts automatically to the possible independence structure of $f$ and this allows us to improve significantly the accuracy of estimation.

\end{abstract}

\end{frontmatter}

\section{Introduction} 
\label{sec:introduction}
Let $X_k=\big(X_{k,1},\ldots,X_{k,d}\big),\;k\in\bN^*,$ be a sequence of $\bR^d$-valued i.i.d. random vectors defined on a complete probability space $\left(\Omega,\mA,\bP\right)$ and having an unknown density $f$ with respect to the Lebesgue measure. Assume that we have at our disposal indirect observations given by
\begin{equation}
\label{eq:model}
Y_k=X_k+\varepsilon_k,\quad k=1,\ldots,n,
\end{equation}
where the errors $\varepsilon_k$ are also i.i.d. $d$-dimensional random vectors, independent of the $X_k$'s, with a known density $q$. 

The goal is to estimate the density $f$ by using observations $Y^{(n)}=(Y_1,\ldots,Y_n)$. By an estimator we mean any $Y^{(n)}$-measurable mapping $\widetilde{f}:\bR^n\rightarrow\bL_p\left(\bR^d\right)$. The accuracy of an estimator is measured by its \textit{$\bL_p$-risk}
$$
\cR_p\left[\widetilde{f},f\right]:=\left(\bE_f\left\|\widetilde{f}-f\right\|_p^p\right)^{\frac{1}{p}},\;\; p\in [1,+\infty),\quad\cR_{\infty}\left[\widetilde{f},f\right]:=\bE_f\left\|\widetilde{f}-f\right\|_{\infty}.
$$
Here and in the sequel $\bE_f$ denotes the expectation with respect to the probability measure $\bP_f$ of the observations $Y^{(n)}=(Y_1,\ldots,Y_n)$ and $\left\|g\right\|_{\mathbf{r}}$ is the $\bL_{\mathbf{r}}$-norm of $g\in\bL_{\mathbf{r}}(\bR^s)$, $s\in\bN^*$, $\mathbf{r}\in [1,+\infty]$. We will also denote by $\widehat{g}$ the Fourier transform of $g\in\bL_1(\bR^s)$, defined by $\widehat{g}(x)=\int e^{i<t,x>}g(x)\rd x$, where $<\cdot ,\cdot >$ is the euclidean scalar product on $\bR^s$.

\smallskip

The aforementioned deconvolution model, which is more realistic than the density model (with direct observations), exists in many different fields and is the subject of many theoretical studies. In most of them, the main interest is to provide estimators which achieve optimal rates of convergence on particular functional classes in a minimax sense. For instance, the problem of minimax estimation in the deconvolution model with pointwise and $\bL_2$ risks was investigated by Carroll and Hall \cite{carroll-hall},  Stefanski \cite{stefanski}, Fan (\cite{Fan1},\cite{Fan2}), Pensky and Vidakovic \cite{pensky-vidakovic}, Butucea \cite{butucea}, Hall and Meister \cite{hall-meister}, Meister \cite{meister}, Butucea and Tsybakov (\cite{tsybakov-butucea1},\cite{tsybakov-butucea2}), Butucea and Comte \cite{comte-butucea}. Global density deconvolution was also considered under a weighted $\bL_p$-norm (defined with an integrable weight function) by Fan \cite{Fan2} and under the sup-norm loss by Stefanski \cite{stefanski}, Bissantz, D\"umgen, Holzmann and Munk \cite{bissantz-al} and Lounici and Nickl \cite{lounicinickl:sup-norm-deconvolution}. Whereas all the works cited above are in the unidimensional setting, the problem of deconvolving a multidimensional density under pointwise or $\bL_2$ loss has been addressed by Masry (\cite{masry1}, \cite{masry2}), Youndj\'e and Wells \cite{youndje-wells} and Comte and Lacour \cite{comte:deconvolution}. 

In the present paper the aim is twofold. First, we deal with optimal deconvolution of a multivariate density under $\bL_p$ and sup-norm losses. Next, as in Lepski \cite{lepski:supnormlossdensityestimation} (under sup-norm loss) and in Rebelles \cite{rebelles2} (under $\bL_p$-losses) for the density model, we also take advantage of the fact that some coordinates of the $X_k$'s may be independent from the others, but in a unified way.

\paragraph{Minimax estimation} In the framework of the minimax estimation it is assumed that $f$ belongs to a certain set  of functions $\Sigma$, and then the accuracy of an estimator $\widetilde{f}$ is measured by its \textit{maximal risk} over $\Sigma$ :
$$
\cR_p\left[\widetilde{f},\Sigma\right]:=\sup_{f\in\Sigma}\cR_p\left[\widetilde{f},f\right].
$$
The objective here is to construct an estimator $\widetilde{f}_*$ which achieves the asymptotic of \textit{the minimax risk} (minimax rate of convergence) :
$$
\cR_p\left[\widetilde{f}_*,\Sigma\right]\asymp\inf_{\widetilde{f}_n}\cR_p\left[\widetilde{f},\Sigma\right]:=\varphi_{n,p}(\Sigma),\quad n\rightarrow +\infty,
$$
where infimum is taken over all possible estimators. Such an estimator is called minimax on $\Sigma$. 

In this paper, we focus on the problem of minimax estimation over anisotropic Nikolskii classes of densities $N_{r,d}(\beta,L)$, see the definition in Section \ref{sec:nikolskii}. Whereas the vector $\beta=(\beta_1,\ldots,\beta_d)$ represents the smoothness of the target density, $r=(r_1,\ldots,r_d)$ represents the index of homogeneity. For the case where $p$ is finite we will assume that the smoothness of $f$ is measured in the same $\bL_p$-norm that the accuracy of estimation, that means $r_j=p$ for $j=1,\ldots,d$. In the latter case, the vector $r$ will be replaced by $p$ in the notation of the functional class. If $\beta_j=\beta_0$, $L_j=L_0$ and $r_j=r_0$ for all $j=1,\ldots,d$, any function belonging to $N_{r_0,d}(\beta_0,L_0)$ is called isotropic function.

\smallskip

In Comte and Lacour \cite{comte:deconvolution} it was shown that
\begin{equation}
\label{eq:L2-rate}
\varphi_{n,2}(N_{2,d}(\beta,L))\asymp n^{-\frac{\tau}{2\tau+1}},\quad\tau:=\left[\sum_{j=1}^d\frac{2\lambda_j+1}{\beta_j}\right]^{-1},
\end{equation}
when the common density $q$ of the errors (which is assumed to be known) satisfies 
$$
\mathbf{A}_1\prod_{j=1}^d\left(1+t^2_j\right)^{-\frac{\lambda_j}{2}}\leq\left|\widehat{q}(t)\right|\leq\mathbf{A}_2\prod_{j=1}^d\left(1+t^2_j\right)^{-\frac{\lambda_j}{2}},\quad\forall t\in\bR^d,
$$
for some constants $\mathbf{A}_1,\mathbf{A}_2,\lambda_j>0$, $j=1,\ldots,d$. Such a density is usually called ordinary smooth of order $\lambda=(\lambda_1,\ldots,\lambda_d)$.

Note that the latter result was proved in the one dimensional setting by Fan \cite{Fan2}. However, whereas Fan \cite{Fan2} provided an estimator whose construction depends on the smoothness parameter $\beta$ of the functional class $N_{2,1}(\beta,L)$ (which is not known in practice), Comte and Lacour \cite{comte:deconvolution} proposed an adaptive strategy. Indeed, they have constructed a single estimator which is fully data driven and minimax on each class $N_{2,d}(\beta,L)$, whatever the nuisance parameter $(\beta,L)$ in a large range. Such an estimator is called optimal adaptive over the scale $\{N_{2,d}(\beta,L)\}_{(\beta,L)}$.

\smallskip

Lounici and Nickl \cite{lounicinickl:sup-norm-deconvolution} considered the problem of adaptive deconvolution of a univariate density under sup-norm loss and proved that
\begin{equation}
\label{eq:sup-norm-rate1}
\varphi_{n,\infty}(N_{\infty,1}(\beta,L))\asymp\left(\frac{n}{\ln(n)}\right)^{-\frac{\tau}{2\tau+1}},\quad\tau:=\left[\frac{2\lambda+1}{\beta}\right]^{-1},
\end{equation}
when the common density $q$ of the errors is ordinary smooth of order $\lambda>0$. Moreover, they provided an optimal adaptive estimator over the scale of H\"older classes $\{N_{\infty,1}(\beta,L)\}_{(\beta,L)}$.

It is worth mentioning that Fan \cite{Fan2}, Lounici and Nickl \cite{lounicinickl:sup-norm-deconvolution} and Comte and Lacour \cite{comte:deconvolution}, as in most of the aforementioned papers, considered also the case of errors having a common density whose Fourier transform has exponential decay, usually called super smooth. In the multidimensional setting, Comte and Lacour \cite{comte:deconvolution} showed that, in presence of super smooth noise, the rates of convergence on anisotropic Nikolskii classes (considered as classes of ordinary smooth densities) are logarithmic and achieved by a kernel estimator whose bandwidth depends only on the smoothness parameters of the noise. Thus, in the latter case, no bandwidth selection procedure is required to get adaptive properties. Note that Youndj\'e and Wells \cite{youndje-wells} considered the problem of adaptive deconvolution of an isotropic density in the ordinary smooth case, namely the "moderately ill-posed" case in inverse problems. The results obtained in Comte and Lacour \cite{comte:deconvolution} under $\bL_2$-loss generalizes considerably those of Youndj\'e and Wells \cite{youndje-wells}.

\smallskip 

In the present paper, we deal with the problem of minimax adaptive deconvolution of an anisotropic density in the ordinary smooth case with $\bL_p$-risks, $p\in[1,\infty]$. The  rates of convergence given in (\ref{eq:L2-rate})-(\ref{eq:sup-norm-rate1}) are recovered from the results we obtain. Indeed, we provide adaptive kernel estimators which achieve the following minimax rates of convergence respectively:
\begin{eqnarray}
\label{eq:Lp-rate}
&&\varphi_{n,p}(N_{p,d}(\beta,L))\asymp n^{-\frac{\tau}{2\tau+1}},\quad\forall p\in[2,+\infty);
\\*[2mm]
\label{eq:sup-rate}
&&\varphi_{n,\infty}(N_{r,d}(\beta,L))\asymp\left(\frac{n}{\ln(n)}\right)^{-\frac{\Upsilon}{2\Upsilon+1}},\quad\Upsilon^{-1}:=\tau^{-1}+[\omega\kappa]^{-1},
\end{eqnarray}
where $\tau$ is given in (\ref{eq:L2-rate}), $\omega:=\left[\sum_{j=1}^d\frac{2\lambda_j+1}{\beta_jr_j}\right]^{-1}$ and $\kappa:=\left(1-\sum_{j=1}^d\frac{1}{\beta_j r_j}\right)\left[\sum_{j=1}^d\frac{1}{\beta_j}\right]^{-1}>0$.

\smallskip

Here, the optimality is a direct consequence of minimax lower bounds recently obtained by Lepski and Willer \cite{lepski:lowerbound-deconvolution}. As usually, these lower bounds hold under additional assumptions on the common density of the errors, see Section \ref{sec:lower-bounds}. Moreover, they proved that there is no uniformly consistent estimator on $N_{r,d}(\beta,L)$ under sup-norm loss if $\kappa\leq 0$. Note also that, for $p\in(1,2)$, our estimator does not achieve the minimax lower bound on $N_{p,d}(\beta,L)$ they have found. Finally, we will not consider the case $p=1$ since the results in Lepski and Willer \cite{lepski:lowerbound-deconvolution} show that there is no uniformly consistent estimator on $N_{1,d}(\beta,L)$ under $\bL_1$-loss.

\smallskip

It is important to emphasize that minimax rates depend heavily on the dimension $d$. To reduce the influence of the dimension on the accuracy of estimation (curse of dimensionality), many researchers have studied the possibility of taking into account, not only the smoothness properties of the target function, but also some structural hypothesis on the statistical model. For instance, see the works on the composite function structure in Horowitz and Mamen \cite{horowitz-mamen}, Iouditski et al. \cite{iouditski-lepski-tsybakov} and Baraud and Birg\'e\cite{baraud}, the works on multi-index structure in Goldenshluger and Lepski \cite{G-L:structural} and Lepski and Serdyukova \cite{lepski:singleindexmodel}, and the works on the multiple index model in density estimation in Samarov and Tsybakov \cite{samarov-tsybakov}.

Below, we discuss one of the possibilities of facing to this problem in the framework of density estimation. The approach which has been recently proposed in Lepski \cite{lepski:supnormlossdensityestimation} is to take into account the independence structure of the target density $f$, namely its product structure due to the independence structure of the vector $X_1$.

\paragraph{Organization of the paper} In Section \ref{sec:assumptions}, we describe assumptions on the densities involved in the statistical model (\ref{eq:model}) and we recall the minimax lower bounds obtained in Lepski and Willer \cite{lepski:lowerbound-deconvolution} useful in this paper. In Section \ref{kernelestimators}, we introduce the family of kernel estimators we use for our procedure and then we describe the selection rule that leads to the construction of our final estimator. In Section \ref{sec:results}, we provide some oracle inequalities and, as consequences, minimax adaptive upper bounds under $\bL_p$-losses over scales of anisotropic Nikolskii classes. Further, we discuss the optimality of our estimator and the influence of the independence structure of the target density on the accuracy of estimation. Proofs of all main results are given in Section \ref{proofs}. Proofs of technical useful results are deferred to the Appendix.

\section{Assumptions on densities $f$ and $q$}
\label{sec:assumptions}

\subsection{Structural assumption on the target density} 
\label{sec:structure}
Denote by $\cI_d$ the set of all subsets of $\{1,\ldots,d\}$, except the empty set. Let $\mP$ be a given set of partitions of $\{1,\ldots,d\}$. For all $I\in \cI_d$ denote also $\overline{I}=\{1,\ldots,d\} \backslash I$ and $\left|I\right|=$card$(I)$. We will use $\overline{\emptyset}$ for $\{1,\ldots,d\}$. Finally, for all $x\in\bR^d$ and $I\in \cI_d$ put $x_I:=(x_i)_{j\in I}$ and, for any probability density $g:\bR^d\rightarrow\bR_+$,
$$
g_I(x_I):=\int_{\bR^{\left|\overline{I}\right|}}g(x)dx_{\overline{I}}.
$$ 
Assume that $g_{\overline{\emptyset}}\equiv g$ and that $g_{\emptyset}\equiv 1$. Note also that $f_I$ and $q_I$ are the marginal densities of $X_{1,I}$ and $\varepsilon_{1,I}$ respectively. 

If $\cP\in\mP$ is such that the vectors $X_{1,I}$, $I\in\cP$, are independent then $f(x)=\prod_{I\in\cP}f_I(x_I),\; \forall x\in\bR^d$. In the sequel, the possible independence structure of the density $f$ will be represented by a partition belonging to the following set :
\begin{equation}
\label{eq:structure}
\mP(f):=\left\{\cP\in\mP:\; f(x)=\prod_{I\in\cP}f_I(x_I),\; \forall x\in\bR^d\right\}.
\end{equation}
Remark that $\mP(f)$ is not empty if we consider that $\overline{\emptyset}\in\mP$, or that $\mP=\{\cP\}$ if the independence structure of $f$ is known. The possibility of choosing $\mP$, instead of considering all partitions of $\{1,\ldots,d\}$, is introduced for technical purposes. This is explained in more detail in Lepski \cite{lepski:supnormlossdensityestimation}, section 2.1, paragraph "\textit{Extra parameters}".

Finally, we endow the set $\mP$ with the operation $"\diamond"$ introduced in Lepski \cite{lepski:supnormlossdensityestimation} : for any $\cP,\cP'\in\mP$ 
\begin{equation}
\label{eq:partionoperation}
\cP\diamond\cP':=\left\{I\cap I'\neq\emptyset,\;I\in\cP,\;I'\in\cP'\right\}.
\end{equation}
The use of this operation for the estimation procedure allows us to construct an estimator which adapts automatically to the independence structure of the underlying density.

\subsection{Noise assumptions for upper bounds} 
\label{sec:noise}
Both the definition of our estimation procedure and the computation of the $\bL_p$-risk, $p\in(1,+\infty]$, lead us to consider that the density $q$ of the noise random vector $\varepsilon_1$ satisfies following assumptions.

\begin{assmN1} Assume that, for any $I\in\cP\diamond\cP'$, $(\cP,\cP')\in\mP\times\mP$:

\smallskip

$(i)$ if $p=2$, then $\left\|\widehat{q_I}\right\|_{1}<+\infty$;

\smallskip

$(ii)$ if $p\in(2,+\infty]$, then $\left\|q_I\right\|_{\infty}<+\infty$.
\end{assmN1}

\smallskip

\begin{assmN2} 
\label{as:errorminimax} Assume that, for some constants $\mathbf{A}\;>0$, $\lambda_j>0$, $j=1,\ldots,d$, one has for any $I\in\cP\diamond\cP'$, $(\cP,\cP')\in\mP\times\mP$:

\smallskip

\noindent $(i)$ if $p=2$,
\begin{eqnarray*}
\label{eq:noiseassumption1}
\left|\widehat{q_I}(t)\right|\geq\mathbf{A}^{-1}\prod_{j\in I}\left(1+t^2_j\right)^{-\frac{\lambda_j}{2}},\quad\forall t\in\bR^d;\nonumber
\end{eqnarray*}

\noindent$(ii)$ if $p\in(1,+\infty)\backslash\{2\}$, $\widehat{q_I}(t_I)\neq 0$, $\forall t\in\bR^d$, $\widehat{q_I}^{\;-1}\in\cC^{|I|}\left(\bR^{|I|}\right)$ and
\begin{eqnarray*}
\label{eq:noiseassumption2}
\left|\left[D^{\alpha_I}\widehat{q_I}^{\;-1}\right](t_I)\prod_{j\in I} t_j^{\alpha_j}\right|\leq\mathbf{A}\prod_{j\in I} \left(1+t^2_j\right)^{\frac{\lambda_j}{2}},\;\forall t\in\bR^d,\quad\forall\alpha_I=(\alpha_j)_{j\in I}\in\bN^{|I|},\;\sum_{j\in I}\alpha_j\leq |I|;\nonumber
\end{eqnarray*}

\noindent$(iii)$ if $p=+\infty$, $\widehat{q_I}(t_I)\neq 0$, $\forall t\in\bR^d$, $\widehat{q_I}^{\;-1}\in\cC^{1}\left(\bR^{|I|}\right)$ and
\begin{eqnarray*}
\label{eq:noiseassumption3}
\left|\left[D_k^{\alpha_k}\widehat{q_I}^{\;-1}\right](t_I)\right|\leq\mathbf{A}\prod_{j\in I}\left(1+t^2_j\right)^{\frac{\lambda_j}{2}},\;\forall t\in\bR^d,\quad\forall k\in I, \forall\alpha_k\in\{0,1\}.\nonumber
\end{eqnarray*}
\end{assmN2}
\noindent Here and in the sequel, $D_k^{\alpha_k} g$ denotes the $\alpha_k$th order partial derivate of $g$ with respect to the $k$th variable, $D_k^{0} g\equiv g$ and, for any multi-index $\alpha=(\alpha_1,\ldots,\alpha_s)\in\bN^s$, $D^{\alpha}g$ denotes the derivative $D_1^{\alpha_1}\ldots D_{s}^{\alpha_{s}}g$ of $g:\bR^s\rightarrow\bR$.

\smallskip

Assumption (N1) is satisfied for many distributions like centered Gaussian, Cauchy, Laplace or Gamma type multivariate ones. Assumption (N2) is quite restrictive since it does not hold for the classical Cauchy and Gaussian densities, whose characteristic functions have exponential decay. However, it is verified by the centered Laplace and Gamma type distributions, whose characteristic functions have polynomial decay. As mentioned in Comte and Lacour \cite{comte:deconvolution}, the latter case keep a great interest in particular physical contexts; see, for instance, the study of the pile-up model in Comte and Rebafka \cite{comte:fluorescencemodel}.

In what follows, we assume that $q$ satisfies Assumptions (N1)-(N2).

\subsection{Smoothness assumption on the target density}
\label{sec:nikolskii}
In the literature there are several definitions of the \textit{anisotropic Nikolskii class of densities} which are equivalent. Let us recall the definition we use in the present paper. Set $\left\{e_1,\ldots,e_s\right\}$, the canonical basis in $\bR^s,\;s\in\bN^*$.

\begin{definition} 
Assume that $r=(r_1,\ldots,r_s)\in[1,+\infty]^s$, $\beta=(\beta_1,\ldots,\beta_s)\in(0,+\infty)^s$ and $L= (L_1,\ldots,L_s)$ $\in(0,+\infty)^s$. A probability density $g:\bR^s\rightarrow\bR_+$ belongs to the anisotropic Nikolskii class $N_{r,s}(\beta,L)$ if 
\begin{eqnarray*}
&&\ (i)\;\left\|D_j^{k} g\right\|_{r_j}\leq L_j,\quad\forall k=0,\ldots,\left\lfloor \beta_j\right\rfloor,\;\;\forall j=1,\ldots,s;
\\*[2mm]
&&\ (ii)\;\left\|D_j^{\left\lfloor \beta_j\right\rfloor} g(\cdot+ze_j)-D_j^{\left\lfloor \beta_j\right\rfloor} g(\cdot)\right\|_{r_j}\leq L_j\left|z\right|^{\beta_j-\left\lfloor \beta_j\right\rfloor},\quad\forall z\in\bR,\;\;\forall j=1,\ldots,s.
\end{eqnarray*}
\end{definition}
Here and in the sequel, $\left\lfloor a\right\rfloor$ is the largest integer strictly less than the real number $a$. Furthermore, we use the notation $N_{\mathbf{r},s}(\beta,L)$ for $N_{r,s}(\beta,L)$ when $r=(\mathbf{r},\ldots,\mathbf{r})$.

\smallskip

In order to take into account the smoothness of the underlying density and its possible independence structure simultaneously, a certain collection of anisotropic Nikolskii classes of densities was introduced in Lepski \cite{lepski:supnormlossdensityestimation}, Section 3, Definition 2. However, since the adaptation is not necessarily considered with respect to the set of all partitions of $\{1,\ldots,d\}$, the condition imposed therein can be weakened. For instance, if $\mP=\{\overline{\emptyset}\}$ (no independence structure), we want to find again the well known results concerning the adaptive estimation over the scale of anisotropic Nikolskii classes of densities $\big\{N_{r,d}(\beta,L)\big\}$, that is not possible with the classes introduced in Lepski \cite{lepski:supnormlossdensityestimation}.  For these reasons, the following collection $\big\{N_{r,d}(\beta,L,\cP)\big\}_{\cP}$ was introduced in Rebelles \cite{rebelles1}, Section 3.1.

\begin{definition}
\label{def:nikolskiidensityadaptive}
Let $r\in[1,+\infty]^d$ and $\left(\beta,L,\cP\right)\in\left(0,+\infty\right)^d\times\left(0,+\infty\right)^d\times\mP$ be fixed. A probability density $g:\bR^d\rightarrow\bR_+$ belongs to the class $N_{r,d}\left(\beta,L,\cP\right)$ if
\begin{eqnarray}
\label{eq:marginalenikolskii2}
g(x)=\prod_{I\in\cP}g_I(x_I),\;\;\forall x\in\bR^d;\quad g_I\in N_{r_I,\left|I\right|}(\beta_I,L_I),\;\;\forall I\in\cP'\diamond\cP'',\;\forall\left(\cP',\;\cP''\right)\in\mP\times\mP.
\end{eqnarray}
\end{definition}

Note that, if $\mP=\{\overline{\emptyset}\}$, the class $N_{r,d}\left(\beta,L,\overline{\emptyset}\right)$ coincides with the classical anisotropic Nikolskii class of densities $N_{r,d}\left(\beta,L\right)$.

\subsection{Noise assumptions for minimax lower bounds} 
\label{sec:lower-bounds}
Recently, Lepski and Willer \cite{lepski:lowerbound-deconvolution} have obtained minimax lower bounds for $\varphi_{n,p}(N_{r,d}(\beta,L))$, $p\in[1,+\infty]$, when the density $q$ of the noise random vector $\varepsilon_1$ satisfies the following assumption. 

\begin{assmN3} 
\label{as:minimaxlowerbound1} For any multi-index $\alpha=(\alpha_1,\ldots,\alpha_d)\in\{0,1\}^{d}$ satisfying $\alpha_1+\ldots+\alpha_d\geq 1$, $D^{\alpha}\widehat{q}$ exists. Furthermore, there exist constants 
$\mathbf{B}\;>0$ and $\lambda_j>0$, $j=1,\ldots,d$, such that: 
\begin{eqnarray*}
& (i)&\left|\widehat{q}(t)\right|\leq\mathbf{B}\prod_{j=1}^d\left(1+t^2_j\right)^{-\frac{\lambda_j}{2}},\quad\forall t\in\bR^d;
\\*[2mm]
& (ii)&\left\|\widehat{q}^{\;-1}D^{\alpha}\widehat{q}\right\|_{\infty}\leq\mathbf{B},\quad\forall\alpha=(\alpha_1,\ldots,\alpha_d)\in\{0,1\}^{d},\;\alpha_1+\ldots+\alpha_d\geq 1.
\end{eqnarray*}
\end{assmN3}

Note first that Assumption (N3) is also verified for centered Laplace or Gamma-type distributions. Next, if $\mP=\{\overline{\emptyset}\}$ (no independence structure), any density $q$ that satisfies both the condition $(i)$ of Assumption (N3) and Assumptions (N2) verifies
$$
\mathbf{A}^{-1}\prod_{j=1}^d\left(1+t^2_j\right)^{-\frac{\lambda_j}{2}}\leq\left|\widehat{q}(t)\right|\leq\mathbf{B}\prod_{j=1}^d\left(1+t^2_j\right)^{-\frac{\lambda_j}{2}},\quad\forall t\in\bR^d,
$$
and hence is ordinary smooth of order $\lambda=(\lambda_1,\ldots,\lambda_2)$. Furthermore, the condition imposed in the left hand side of the latter inequalities, together with the condition $(ii)$ of Assumption (N3) (or Condition 1 in Lounici and Nickl \cite{lounicinickl:sup-norm-deconvolution} for the one dimensional setting), implies that condition $(iii)$ of Assumption (N2) is satisfied.

The results below follow from Theorems 2 and 3 in Lepski and Willer \cite{lepski:lowerbound-deconvolution} and allow us to assert the optimality of our estimators when $\mP=\{\overline{\emptyset}\}$ (no independence structure).

\begin{theorem} 
\label{theo:minimawlowerbound1}
Let $L_0>0$ and $p\in[2,+\infty)$ be fixed. Suppose that Assumptions (N3) is satisfied. 
Then, for any $\left(\beta,L\right)\in(0,\infty)^d\times[L_0,\infty)^d$
$$
\liminf_{n\rightarrow +\infty}\inf_{\widetilde{f}_n}\left\{\varphi_{n,p}^{-1}(N_{p,d}(\beta,L))\cR_p\left[\widetilde{f}_n,N_{p,d}\left(\beta,L\right)\right]\right\}>0,
$$
where infimum is taken over all possible estimators and $\varphi_{n,p}(N_{p,d}(\beta,L))$ is given in (\ref{eq:Lp-rate}).
\end{theorem}

\begin{theorem} 
\label{theo:minimawlowerbound2}
Let $L_0>0$ and $\left(\beta,L,r\right)\in(0,\infty)^d\times[L_0,\infty)^d\times[1,\infty]^d$ be fixed. Suppose that Assumptions (N3) is satisfied. Then,

\smallskip
 
\noindent $(i)$ there is no uniformly consistent estimator if $1-\sum_{j=1}^d\frac{1}{\beta_j r_j}\leq 0$;

\smallskip

\noindent $(ii)$ if $1-\sum_{j=1}^d\frac{1}{\beta_j r_j}>0$
$$
\liminf_{n\rightarrow +\infty}\inf_{\widetilde{f}_n}\left\{\varphi_{n,\infty}^{-1}(N_{r,d}(\beta,L))\cR_{\infty}\left[\widetilde{f}_n,N_{r,d}\left(\beta,L\right)\right]\right\}>0,
$$
where infimum is taken over all possible estimators and $\varphi_{n,\infty}(N_{r,d}(\beta,L))$ is given in (\ref{eq:sup-rate}).
\end{theorem}

\section{Estimation procedure}
\label{kernelestimators}

In this section, we construct an estimator following a scheme of selection rule introduced in Lepski \cite{lepski:supnormlossdensityestimation} to take into account the possible independence structure of the underlying density. If $\mP=\{\overline{\emptyset}\}$ this scheme coincides with a version of the methodology proposed by Goldenshluger and Lepski \cite{G-L:density-Lploss}. This methodology, employed in many areas of nonparametric statistics, has been recently used by Comte and Lacour \cite{comte:deconvolution} in the framework of the deconvolution model.

\subsection{Kernel-type estimators}
\label{sec:kernelestimators}

Let $\textbf{K}:\bR\rightarrow\bR$ be a fixed symmetric kernel ($\int{\textbf{K}}=1$) belonging to the well known Schwartz class $\bS(\bR)$. For instance, $\textbf{K}$ may be a Gaussian kernel. For all $I\in \cI_d$, $h\in (0,1]^d$ and $x\in\bR^d$ put
$$
K_I(x_I):=\prod_{j\in I}\textbf{K}(x_j),\quad K_{h_I}(x_I):=V_{h_I}^{-1}\prod_{j\in I}\textbf{K}(x_j/h_j),\quad V_{h_I}:=\prod_{j\in I}h_j. 
$$
Therefore, in view of the definition of both the kernel $\textbf{K}$ and Assumption (N2) on the errors, one can define the \textit{kernel-type estimator}
\begin{eqnarray}
\label{eq:kernelestimator}
\widetilde{f}_{h_I}(x_I):=n^{-1}\sum_{k=1}^{n}L_{(h_I)}\left(x_I-Y_{k,I}\right),\quad L_{(h_I)}(x_I):=\frac{1}{(2\pi)^{|I|}}\int_{\bR^{|I|}}e^{-i\left\langle t_I,x_I\right\rangle}\frac{\widehat{K_{h_I}}(t_I)}{\widehat{q_I}(t_I)}\rd t_I.
\end{eqnarray} 

\smallskip

The ideas that led to the introduction of the estimators $\widetilde{f}_{h_I}$ are explained in Fan \cite{Fan2} in the one-dimensional setting and, in Comte and Lacour \cite{comte:deconvolution} in the multivariate context.

\smallskip

\paragraph{\textbf{Family of estimators}} Below we propose a data driven selection from the family of estimators
\begin{equation}
\label{eq:estimatorsfamily}
\mF\left[\;\mP\;\right]:=\left\{\widetilde{f}_{(h,\cP)}(x)=\prod_{I\in\cP}\widetilde{f}_{h_I}(x_{I}),\;x\in\bR^d,\;(h,\cP)\in\cH_p[\;\mP\;]\right\},
\end{equation}
where the set $\cH_p[\;\mP\;]$ of parameters $(h,\cP)$ is constructed as follows.

\smallskip 

For $I\in\cI_d$, consider first the set of multibandwidths 
$$
\mH_{p,I}:=\left\{h_I\in\left[h_{min}^{(p)},h_{max}^{(p)}\right]^{|I|}:\; h_j=2^{-k_j},\;k_j\in\bN^*,\;j\in I\right\},
$$
\begin{eqnarray*}
h_{min}^{(p)}:=\left\{
\begin{array}{ll}
n^{-\left(1\vee\frac{p}{|I|}\right)}, &\ p\in(1,+\infty),
\\*[2mm]
n^{-1}, &\ p=+\infty,
\end{array}
\right.
\quad h_{max}^{(p)}:=\left\{
\begin{array}{ll}
[\ln(n)]^{-\frac{p}{|I|}}, &\ p\in(1,+\infty),
\\*[2mm]
1, &\ p=+\infty.
\end{array}
\right.
\end{eqnarray*}
Then define
\begin{equation}
\label{eq:multibandwidthset}
\cH_{p,I}:=\left\{h_I\in\mH_{p,I}:\; \left(nV_{h_I}\right)^{b_p}\prod_{j\in I}h_j^{\lambda_j}\geq c_{p}\textbf{1}_{\{p<\infty\}}+\sqrt{\ln(n)}\textbf{1}_{\{p=+\infty\}}\right\},
\end{equation}
$$
b_p:=\frac{1}{2}\wedge(1-\frac{1}{p}),\;\;c_{p}:=1\wedge \left\{\frac{p}{e}\left[1+\lambda_{max}\left(2\vee\frac{p}{p-1}\right)\right]\right\}^{-\;p\left[b_p+\lambda_{max}\right]},\;\;\lambda_{max}:=\max_{j=1,\ldots,d}\lambda_{j}.
$$
The constant $c_p$ is chosen in order to have $\cH_{p,I}\neq\emptyset,\;\forall n\geq 3$. 

\smallskip

Put finally
$$
\cH_p[\;\mP\;]:=\left\{(h,\cP)\in (0,1]^d\times\mP:\; h_I\in\cH_{p,I},\;\forall I\in\cP\;\right\}.
$$

The introduction of the estimator $\widetilde{f}_{(h,\cP)}$ is based on the following simple observation. If there exists $\cP\in\mP(f)$,  the idea is to estimate separately each marginal density corresponding to $I\in\cP$. Since the estimated density possesses the product structure we seek its estimator in the same form.

\smallskip

\paragraph{\textbf{Auxiliary estimators}} We mimic the procedure of Lepski \cite{lepski:supnormlossdensityestimation} by introducing the following auxiliary estimators. Consider first the classical kernel auxiliary estimators 
$$
\widetilde{f}_{h_I,\eta_I}(x_I):=K_{\eta_I}\star \widetilde{f}_{h_I}(x_I),\;h,\eta\in(0,1]^d,\;I\in\cI_d,
$$
where, here and in the sequel, $"\star"$ stands for the standard convolution product on $\bR^{s}$, $s\in\bN^*$. 

Then put, for $h,\eta\in(0,1]^d$ and $\cP,\cP'\in\mP$,
$$
\widetilde{f}_{(h,\cP),(\eta,\cP')}(x):=\prod_{I\in\cP\diamond\cP'}\widetilde{f}_{h_{I},\eta_{I}}(x_{I}),
$$
where the operation $"\diamond"$ is defined by (\ref{eq:partionoperation}).

The ideas that led to the introduction of the estimators $\widetilde{f}_{(h,\cP),(\eta,\cP')}$, based on both the operation "$\star$" and "$\diamond$", are explained in Lepski \cite{lepski:supnormlossdensityestimation}, Section 2.1, paragraph "\textsl{Estimation construction}". Note that the arguments given in the latter paper do not depend on the norm used in the definition of the risk and remain valid for estimation under $\bL_p-$loss.

\subsection{Selection rule}
\label{sec:selectionrule}

For $I\in\cI_d$ and $h\in(0,1]^d$, define
\begin{eqnarray*}
\cU_{p}(h_I):=\left\{
\begin{array}{llll}
n^{\frac{1}{p}-1}\left\|L_{(h_I)}\right\|_{p}, &\ p\in(1,2),
\\*[4mm]
n^{-\frac{1}{2}}\prod_{j\in I}h_j^{-\lambda_j-\frac{1}{2}}, &\ p=2,
\\*[4mm]
n^{-\frac{1}{2}}\left[\prod_{j\in I}h_j^{-\lambda_j-\frac{1}{2}}+\sqrt{\ln(n)}\left\|L_{(h_I)}\right\|_{\frac{2p}{p+2}}\right], &\ p\in(2,+\infty),
\\*[4mm]
n^{-\frac{1}{2}}\sqrt{\ln(n)}\prod_{j\in I}h_j^{-\lambda_j-\frac{1}{2}}, &\ p=+\infty.
\end{array}
\right.
\end{eqnarray*}

\smallskip 

Put also $\Lambda_{p}:=\md\gamma_p\left[\;\overline{G}_{p}\right]^{\md(\md-1)},$ where $\md:=\sup_{\cP\in\mP}|\cP|$, 
$$
\overline{G}_{p}:=1\vee\left[\left\|\textbf{K}\right\|_{1}^d\sup_{(h,\cP)\in\;\cH_p[\;\mP\;]}\sup_{\cP'\in\;\mP\;}\sup_{I\in\cP\diamond\cP'}\left\|\widetilde{f}_{h_I}\right\|_{p}\right]
$$
and $\gamma_p>0$ is a numerical constant whose expression is given in Section \ref{sec:notations} below.

\smallskip

For $h\in(0,1]^d$ and $\cP\in\mP$ introduce $\cU_{p}(h,\cP):=\sup_{I\in\cP}\cU_{p}(h_I)$ and
\begin{equation}
\label{eq:selectionrule1}
\widetilde{\Delta}_{p}(h,\cP):=\sup_{(\eta,\cP')\in\cH_p\left[\;\mP\;\right]}\left[\left\|\widetilde{f}_{(h,\cP),(\eta,\cP')}-\widetilde{f}_{(\eta,\cP')}\right\|_p-\Lambda_{p}\cU_{p}(\eta,\cP')\right]_+.
\end{equation}
Define finally $\big(\widetilde{h},\widetilde{\cP}\big)$ satisfying
\begin{equation}
\label{eq:selectionrule2}
\widetilde{\Delta}_{p}(\widetilde{h},\widetilde{\cP})+\Lambda_{p}\cU_{p}(\widetilde{h},\widetilde{\cP})=\inf_{(h,\cP)\in\cH_p\left[\;\mP\;\right]}\left[\widetilde{\Delta}_{p}(h,\cP)+\Lambda_{p}\cU_{p}(h,\cP)\right].
\end{equation}
Our selected estimator is $\widetilde{f}:=\widetilde{f}_{(\widetilde{h},\widetilde{\cP})}$.

\smallskip

Note first that the existence of the quantities involved in the selection procedure is ensured by both the finiteness of the set $\cH_p\left[\;\mP\;\right]$ and the following result. The first statement given in Proposition \ref{multipliertheorem} is a simple consequence of Marcinkiewicz Multiplier Theorem; see Theorem 5.2.4. and Corollary 5.2.5. in Grafakos \cite{Grafakos:multiplier}.

\begin{proposition}
\label{multipliertheorem}
Assume that Assumptions (N1)-(N2) are satisfied.

\smallskip

(i) For any $p\in (1,2)$ and any $I\in\cP\diamond\cP'$, $(\cP,\cP')\in\mP\times\mP$, there exists a constant $C_{p,I}:=C_{p,I}(|I|,\mathbf{K},q)>0$
$$
\left\|L_{h_I}\right\|_{p}\leq C_{p,I}\left(V_{h_I}\right)^{-(1-1/p)}\prod_{j\in I}h_j^{-\lambda_j},\quad\forall h\in (0,1]^d.
$$

\smallskip

(ii) For any $I\in\cP\diamond\cP'$, $(\cP,\cP')\in\mP\times\mP$, there exists a constant $C_I:=C_I(|I|,\mathbf{K},q)>0$ such that
\begin{eqnarray*}
\left\|L_{h_I}\right\|_{2}\leq C_I\prod_{j\in I} h_j^{-\lambda_j-\frac{1}{2}},\quad\left\|L_{h_I}\right\|_{\infty}\leq C_I\prod_{j\in I} h_j^{-\lambda_j-1},\; \forall h\in (0,1]^d.
\end{eqnarray*}
\end{proposition}
\noindent The proof of this proposition is postponed to Appendix. It is important to emphasize that the first bound was not used for the definition of $\cU_p(h_I)$ since a dimensional constant is not explicitly done in Theorem 5.2.4. of Grafakos \cite{Grafakos:multiplier}.

Next, we also emphasize that the quantity $\cU_{p}(h_I)$ can be viewed, up to a numerical constant, as a uniform bound on the $\bL_p-$norm of the stochastic error provided by the kernel-type estimator $\widetilde{f}_{h_I}$. This is explained by the following result. For $I\in\cI_d$, $h\in(0,1]^d$ and $x\in\bR^d$, define
$$
\xi_{h_I}(x_I):=\widetilde{f}_{h_I}(x_I)-\bE\{\widetilde{f}_{h_I}(x_I)\}.
$$

\begin{proposition}
\label{prop:empiricalupperbound1} Assume that Assumptions (N1)-(N2) are verified. Let $I\in\cP\diamond\cP'$, $(\cP,\cP')\in\mP\times\mP$, be arbitrary fixed.
If $p\in(1,+\infty]$, $\mathbf{r}\geq 1$ and $n\geq 3$ then
\begin{eqnarray}
\left\{\bE\sup_{h_I\in\cH_{p,I}}\left[\left\|\xi_{h_I}\right\|_{p}-\gamma_{p,I}(\mathbf{r})\cU_{p}(h_I)\;\right]_+^{\mathbf{r}}\right\}^{\frac{1}{\mathbf{r}}}
\leq c_p(\mathbf{r}) n^{-\frac{1}{2}},\quad c_p(\mathbf{r})>0.
\end{eqnarray}
\end{proposition}
\noindent The constants $\gamma_{p,I}(\mathbf{r})$ and $c_p(\mathbf{r})$ do not depend on the sample size $n$. Their explicit expressions can be found in the proof of the latter result, which is also postponed to Appendix.

\smallskip

Finally, in view of the assumptions on the kernel $\textbf{K}$, since $\cH_p\big[\;\mP\;\big]$ is a finite set, $\big(\widetilde{h},\widetilde{\cP}\big)$ exists, is in $\cH_p\big[\;\mP\;\big]$ and is $Y^{(n)}-$measurable. It follows that $\widetilde{f}:\bR^n\rightarrow\bL_p\left(\bR^d\right)$ is an $Y^{(n)}-$measurable mapping.

\section{Main results}
\label{sec:results}
In this section, we first provide oracle inequalities for our estimator $\widetilde{f}$. Then, we discuss adaptive minimax estimation over scales of anisotropic Nikolskii classes.

\subsection{Oracle inequalities}
\label{sec:oracleinequality}

Note that the construction of the proposed procedure does not require any condition concerning the density $f$. However, the following mild assumption will be used for computing its risk: 
\begin{equation}
\label{eq:densityassumption}
f\in\bF_p\left[\;\mP\;\right]:=\left\{g\in\bF:\;\sup_{\cP ,\cP'\in\mP}\sup_{I\in\cP\diamond\cP'}\left\|g_I\right\|_{p}<\infty\;\right\},
\end{equation}
where $\bF$ denotes the set of all probability densities $g:\bR^d\rightarrow\bR_+$. The considered class of densities is determined by the choice of $\mP$ and in particular
$$
\bF_p\left[\;\big\{\overline{\emptyset}\big\}\;\right]=\Big\{g\in\bF:\;\left\|g\right\|_{p}<\infty\Big\},\quad\bF_p\left[\;\left\{\cP\right\}\;\right]=\Big\{g\in\bF:\;\sup_{I\in\cP}\left\|g_I\right\|_{p}<\infty\;\Big\}.
$$

Define, for $(h,\cP)\in\cH_p\left[\;\mP\;\right]$ such that $\cP\in\mP(f)$,
\begin{eqnarray*}
&&\cR_{p}\left[(h,\cP),f\right]:=\left(\bE_f\sup_{\cP'\in\mP}\sup_{I\in\cP\diamond\cP'}\left\|\widetilde{f}_{h_{I}}-f_{I}\right\|_{p}^p\right)^{\frac{1}{p}},\; p\in(1,+\infty),
\\*[2mm]
&&\cR_{\infty}\left[(h,\cP),f\right]:=\bE_f\sup_{\cP'\in\mP}\sup_{I\in\cP\diamond\cP'}\left\|\widetilde{f}_{h_{I}}-f_{I}\right\|_{\infty}.
\end{eqnarray*}

If the possible independence structure $\cP$ of the target density is known, the latter quantity can be viewed as an "$\bL_p-$\textit{risk}" of the estimator $\widetilde{f}_{(h,\cP)}$, defined with the loss
$$
l\left(\widetilde{f}_{(h,\cP)},f\right):=\sup_{\cP'\in\mP}\sup_{I\in\cP\diamond\cP'}\left\|\widetilde{f}_{h_{I}}-f_{I}\right\|_{p}.
$$
In this case, we see that the effective dimension of estimation is not $d$, but $d(\cP):=\sup_{I\in\cP}\left|I\right|$. Therefore, the best estimator from the family $\mF[\;\mP\;]$ (the oracle) should be $\widetilde{f}_{(h^*,\cP^*)}$ such that 
$$
\cR_{p}\left[(h^*,\cP^*),f\right]=\inf_{(h,\cP)\in\cH_p\left[\;\mP\;\right]:\cP\in\mP(f)}\cR_{p}\left[(h,\cP),f\right].
$$

Let us provide the following oracle inequalities for our selected estimator $\widetilde{f}$.

\begin{theorem}
\label{theo:oracleinequality1} Assume that Assumptions (N1)-(N2) are satisfied. 

\noindent If $n\geq 3$ and $p\in(1,+\infty]$ then: $\forall f\in\bF_p\left[\mP\right]$,
\begin{eqnarray}
\label{eq:oracleinequality1}
\cR_p\left[\widetilde{f},f\right]\leq \mathbf{C}_{p,1}(\mathbf{f}_p)\inf_{(h,\cP)\in\cH_p\left[\;\mP\;\right]:\cP\in\mP(f)}\left\{\cR_{p}\left[(h,\cP),f\right]+\gamma_p\cU_{p}(h,\cP)\right\}+\mathbf{C}_{p,2}(\mathbf{f}_p)n^{-\frac{1}{2}},
\end{eqnarray}
where $\mathbf{f}_p:=1\vee\left[\sup_{\cP ,\cP'\in\mP}\sup_{I\in\cP\diamond\cP'}\left\|f_I\right\|_{p}\right]$.
\end{theorem}
\noindent The explicit expression of $\mathbf{C}_{p,1}(\mathbf{f}_p)=\mathbf{C}_{p,1}(d,\mP,\mathbf{K},q,\mathbf{f}_p)$ and $\mathbf{C}_{p,2}(\mathbf{f}_p)=\mathbf{C}_{p,2}(d,\mP,\mathbf{K},q,\mathbf{f}_p)$ is given in the proof of the theorem. It is worth to note that the maps $\mathbf{f}_p\mapsto\mathbf{C}_{p,1}(\mathbf{f}_p)$ and $\mathbf{f}_p\mapsto\mathbf{C}_{p,2}(\mathbf{f}_p)$ are bounded on any bounded interval of $\bR_+$.

\smallskip

If $\mP=\big\{\overline{\emptyset}\big\}$ we obtain automatically some oracle inequalities for estimation on $\bR^d$ under $\bL_p-$loss, without considering any independence structure. In this case, the result above can be improved. Indeed, by scrutinizing its proof, one can easily see that the following theorem is true.

\begin{theorem}
\label{theo:oracleinequality2} Assume that $\mP=\big\{\overline{\emptyset}\big\}$ and that Assumptions (N1)-(N2) are satisfied. 

\noindent If $n\geq 3$ and $p\in(1,+\infty]$ then: $\forall f\in\bF$,
\begin{eqnarray}
\label{eq:oracleinequality3}
\cR_p\left[\widetilde{f},f\right]\leq \inf_{h\in\cH_{p,\overline{\emptyset}}}\left\{\left(1+2\left\|\textbf{K}\right\|_{1}^d\right)\cR_p\left[\widetilde{f}_h,f\right]
+2\gamma_p\cU_{p}(h)\right\}+2\mathbf{C}_{p} n^{-\frac{1}{2}}.
\end{eqnarray}
\end{theorem}
\noindent The explicit expression of the absolute constant $\mathbf{C}_{p}=\mathbf{C}_{p}(d,\mP,\mathbf{K},q)>0$ is given in the proof of the theorem.

\smallskip

Note first that the statement of Theorem \ref{theo:oracleinequality2} holds for all probability densities $f\in\bF$, that is not true  for Theorem \ref{theo:oracleinequality1}. Next, the constant $1+2\left\|\textbf{K}\right\|_{1}^d$ is more suitable than $\mathbf{C}_{p,1}(\mathbf{f}_p)$. Indeed, the prime interest in the oracle approach is to obtain a constant that does not depend on the target density and close to one. However, Theorem \ref{theo:oracleinequality1} allows us to consider both the smoothness properties and the independence structure of the target density and then to reduce the influence of the dimension on the accuracy of estimation. Indeed, if $f$ has an independence structure $\cP\neq\overline{\emptyset}$ and the smoothness parameter $h$ is fixed and properly chosen then our procedure should choose the true partition $\cP$ and the estimator $\widetilde{f}_{(h,\cP)}$ should provide a better accuracy of estimation than the classical kernel-type estimator $\widetilde{f}_{h}$. This was illustrated by a short simulation study in Rebelles \cite{rebelles2} for the density model (with direct observations), under the $\bL_2$-loss.

\subsection{$\bL_p$-adaptive minimax estimation}
\label{minimax-adaptive-estimation}

In what follows, we illustrate the application of Theorems \ref{theo:oracleinequality1} and \ref{theo:oracleinequality2} to adaptive estimation over anisotropic Nikolskii classes of densities $N_{r,d}\left(\beta,L,\cP\right)$ and $N_{r,d}\left(\beta,L\right)$ respectively. To compute an $\bL_p$-risk of a kernel-type estimator, we first compute its bias. Thus, we need to enforce the assumptions imposed on the kernel $\textbf{K}$. One of the possibilities is the following, proposed in Kerkyacharian, Lepski and Picard \cite{kerk2}.

For a given integer $l\geq 2$ and a given symmetric function $u:\bR\rightarrow\bR$ belonging to the Schwartz class $\bS(\bR)$ and satisfying $\int_{\bR} u(z)\rd z=1$ set
\begin{equation}
\label{eq:orthogonalkernel}
u_l(z):=\sum_{j=1}^{l}
\left(\begin{array}{cc}
l\\
j
\end{array}\right)
(-1)^{j+1}\frac{1}{j}u\left(\frac{z}{j}\right),\quad z\in\bR.
\end{equation}
Furthermore we use $\textbf{K}\equiv u_l$ in the definition of the collection of estimators $\mF[\mP]$.\\
The relation of kernel $u_l$ to anisotropic Nikolskii classes is discussed in Kerkyacharian, Lepski and Picard \cite{kerk2}. In particular, it has been shown that
\begin{equation}
\label{eq:kernelorthogonality}
\int_{\bR}\textbf{K}(z)\rd z=1,\qquad\int_{\bR}z^k\textbf{K}(z)\rd z=0,\quad\forall k=1,\ldots,l-1.
\end{equation}

\subsubsection{Minimax adaptive estimation under an $\bL_p$-loss} For $\left(\beta,\cP\right)\in\left(0,+\infty\right)^d\times\mP$ define $\phi_{n,p}\left(\beta,\cP\right):=n^{-\frac{\tau}{2\tau+1}}$, where
\begin{gather}
\label{eq:Lpminimaxrate}
\tau:=\tau(\beta,\cP)=\inf_{I\in\cP}\tau_I,\qquad \tau_I:=\left[\sum_{j\in I}\frac{b_p^{-1}\lambda_j+1}{\beta_j}\right]^{-1},
\end{gather}
where $b_p$ is given in (\ref{eq:multibandwidthset}).

Assume that $\overline{\emptyset}\in\mP$ and consider the estimator $\widetilde{f}$ defined by the selection rule (\ref{eq:selectionrule1})-(\ref{eq:selectionrule2}) with $p\in(1,+\infty)$.

\begin{theorem} 
\label{theo:Lpadaptiveupperbound}
Let $p\in(1,+\infty)$ be arbitrary fixed. Suppose that Assumptions (N1)-(N2) are satisfied. Then for any $\left(\beta,L,\cP\right)\in(0,l]^d\times(0,\infty)^d\times\mP$ one has
$$
\limsup_{n\rightarrow +\infty}\left\{\phi_{n,p}^{-1}(\beta,\cP)\cR_p\left[\widetilde{f},N_{p,d}\left(\beta,L,\cP\right)\right]\right\}<\infty.
$$
\end{theorem}

To get the statement of the latter theorem we apply Theorem \ref{theo:oracleinequality1}. If $\mP=\{\overline{\emptyset}\}$ (no independence structure), we obtain the following theorem by applying Theorem \ref{theo:oracleinequality2}.

\begin{theorem} 
\label{theo:Lpadaptiveupperbound2}
Let $p\in(1,+\infty)$ be arbitrary fixed. Suppose that $\mP=\{\overline{\emptyset}\}$ and that Assumptions (N1)-(N2) are satisfied. Then for any $\left(\beta,L\right)\in(0,l]^d\times(0,\infty)^d$ one has
$$
\limsup_{n\rightarrow +\infty}\left\{\phi_{n,p}^{-1}(\beta,\overline{\emptyset})\cR_p\left[\widetilde{f},N_{p,d}\left(\beta,L\right)\right]\right\}<\infty.
$$
\end{theorem}
To the best of our knowledge, the latter results are new. In view of the assertion of Theorem \ref{theo:minimawlowerbound1}, if $p\in[2,+\infty)$ and Assumptions (N1)-(N3) on the errors are satisfied, we deduce from Theorem \ref{theo:Lpadaptiveupperbound2} that $\phi_{n,p}(\beta,\overline{\emptyset})$ is the minimax rate of convergence on the anisotropic Nikolskii class $N_{p,d}\left(\beta,L\right)$ and that a minimax estimator can be selected from the collection of kernel-type estimators introduced in Section \ref{sec:kernelestimators}. Moreover, if $\mP=\{\overline{\emptyset}\}$ (no independence structure), the quality of estimation of our estimator $\widetilde{f}$ is optimal, up to a numerical constant, on each class $N_{p,d}\left(\beta,L\right)$, whatever the nuisance parameter $\left(\beta,L\right)$. Thus, in the aforementioned case, $\widetilde{f}$ is an optimal adaptive estimator over the scale $\{N_{p,d}\left(\beta,L\right)\}_{(\beta,L)}$. 

Remark that $\bL_p$-estimation of an anisotropic density in the deconvolution model does not require that this density is uniformly bounded, whereas it is imposed in all the works concerning the density model (with direct observations); see, e.g. in Goldenshluger and Lepski \cite{G-L:density-Lploss}.

Unfortunately, if $p\in(1,2)$, our estimator does not achieve the minimax lower bound on $N_{p,d}\left(\beta,L\right)$ obtained in Lepski and Willer \cite{lepski:lowerbound-deconvolution} under the $\bL_p$-loss. We conclude that either our estimator is not minimax on $N_{p,d}\left(\beta,L\right)$ or the lower bound in Lepski and Willer \cite{lepski:lowerbound-deconvolution} is not the minimax rate of convergence on the latter functional class.

It is important to emphasize that both Theorems \ref{theo:Lpadaptiveupperbound}-\ref{theo:Lpadaptiveupperbound2} allow us to analyze the influence of the independence structure on the accuracy of estimation under an $\bL_p$-loss in the deconvolution model. Indeed, we see that
$$
\phi_{n,p}(\beta,\overline{\emptyset})\gg\phi_{n,p}(\beta,\cP),\quad\cP\neq\overline{\emptyset},
$$
whatever the independence structure of the common density of the errors. Thus, our estimation procedure allows us to improve significantly the accuracy of estimation if the target density has an independence structure $\cP\neq\overline{\emptyset}$.

Having said that, the question is: is $\phi_{n,p}(\beta,\cP)$ the minimax rate of convergence on the functional class  $N_{p,d}\left(\beta,L,\cP\right)$? For the density model (that corresponds to $\lambda_j=0$, $j=1,\ldots,d$), it is asserted in Rebelles \cite{rebelles2} that the answer is positive and that the proof of the corresponding minimax lower bound coincides with the one of Theorem 3 in Goldenshluger and Lepski \cite{G-L:density-Lploss-bis}, up to minor modifications to take into account the independence structure. For the deconvolution model, we conjecture that the answer is also positive if $p\in[2,+\infty)$ and that a minimax lower bound on $N_{p,d}\left(\beta,L,\cP\right)$ can be obtained, up to straightforward modifications, as in Lepski and Willer \cite{lepski:lowerbound-deconvolution}.

\subsubsection{Minimax adaptive estimation under sup-norm loss}
\label{sec:supminimax}
For $\left(\beta,r,\cP\right)\in\left(0,+\infty\right)^d\times[1,+\infty]^d\times\mP$ define $\phi_{n,\infty}\left(\beta,r,\cP\right):=\left(\frac{n}{\ln(n)}\right)^{-\frac{\Upsilon}{2\Upsilon+1}}$, where
$$
\Upsilon:=\Upsilon(\beta,r,\cP)=\inf_{I\in\cP}\Upsilon_I,\qquad\Upsilon_I:=\left(\tau_I^{-1}+\left[\omega_I\kappa_I\right]^{-1}\right)^{-1},
$$
\begin{eqnarray}
\label{eq:suprate}
\tau_I:=\left[\sum_{j\in I}\frac{2\lambda_j+1}{\beta_j}\right]^{-1},\quad\omega_I:=\left[\sum_{j\in I}\frac{2\lambda_j+1}{\beta_j r_j}\right]^{-1},\quad\kappa_I:=\frac{1-\sum_{j\in I}\frac{1}{\beta_j r_j}}{\sum_{j\in I}\frac{1}{\beta_j}}.
\end{eqnarray}

Assume that $\overline{\emptyset}\in\mP$ and consider the estimator $\widetilde{f}$ defined by the selection rule (\ref{eq:selectionrule1})-(\ref{eq:selectionrule2}) with $p=+\infty$. As previously, we obtain the following two theorems:

\begin{theorem} 
\label{theo:supadaptiveupperbound}
Suppose that Assumptions (N1)-(N2) are satisfied. Then for any $\left(\beta,L,r,\cP\right)\in(0,l]^d\times(0,\infty)^d\times[1,+\infty]^d\times\mP$ satisfying $1-\sum_{j=1}^d\frac{1}{\beta_j r_j}>0$ one has
$$
\limsup_{n\rightarrow +\infty}\left\{\phi_{n,\infty}^{-1}(\beta,r,\cP)\cR_p\left[\widetilde{f},N_{r,d}\left(\beta,L,\cP\right)\right]\right\}<\infty.
$$
\end{theorem}

\begin{theorem} 
\label{theo:supadaptiveupperbound2}
Suppose that $\mP=\{\overline{\emptyset}\}$ and that Assumptions (N1)-(N2) are satisfied. Then for any $\left(\beta,L,r\right)\in(0,l]^d\times(0,\infty)^d\times[1,+\infty]^d$ satisfying $1-\sum_{j=1}^d\frac{1}{\beta_j r_j}>0$ one has
$$
\limsup_{n\rightarrow +\infty}\left\{\phi_{n,\infty}^{-1}(\beta,r,\overline{\emptyset})\cR_p\left[\widetilde{f},N_{r,d}\left(\beta,L\right)\right]\right\}<\infty.
$$
\end{theorem}
To the best of our knowledge, the latter results are also new. Note first that in the case of direct observations we find again the results obtained in Lepski \cite{lepski:supnormlossdensityestimation}. Next, if $\mP=\{\overline{\emptyset}\}$ and $1-\sum_{j=1}^d\frac{1}{\beta_j r_j}>0$, it follows from Theorems \ref{theo:minimawlowerbound2} and \ref{theo:supadaptiveupperbound2} that, in presence of the noise satisfying Assumptions (N1)-(N3), $\phi_{n,\infty}(\beta,r,\overline{\emptyset})$ is the minimax rate of convergence on the anisotropic class $N_{r,d}\left(\beta,L\right)$. In this case, our estimator is an optimal adaptive one over the scale 
$$
\left\{N_{r,d}\left(\beta,L\right),\;\left(\beta,L,r\right)\in(0,l]^d\times(0,\infty)^d\times[1,+\infty]^d,\; 1-\sum_{j=1}^d\frac{1}{\beta_j r_j}>0\right\}.
$$

It is worth to note that our estimator can be used for pointwise estimation. Moreover, it follows from Theorem \ref{theo:supadaptiveupperbound2} that our estimator achieves the adaptive rates of convergence found in Comte and Lacour \cite{comte:deconvolution} with a pointwise criterion over the scale of H\"older classes $\{N_{\infty,d}\left(\beta,L\right)\}_{(\beta,L)}$.

As previously, Theorems \ref{theo:supadaptiveupperbound}-\ref{theo:supadaptiveupperbound2} allow us to conclude that our procedure leads to a better accuracy of estimation under sup-norm loss whenever the target density has an independence structure $\cP\neq\overline{\emptyset}$. In this case, we improve significantly the results obtained in Comte and Lacour \cite{comte:deconvolution} under a pointwise loss. Furthermore, we emphasize that if the target density has a known independence structure $\cP$, $\mP=\{\cP\}$, $1-\sum_{j=1}^d\frac{1}{\beta_j r_j}\leq 0$ and $1-\sum_{j\in I}^d\frac{1}{\beta_j r_j}>0$, $\forall I\in\cP$, our estimator achieves the rate of convergence $\phi_{n,\infty}(\beta,r,\cP)$ on $N_{r,d}\left(\beta,L,\cP\right)$ whereas there is no uniformly consistent estimator on $N_{r,d}\left(\beta,L\right)$.

Finally, we conjecture that $\phi_{n,\infty}(\beta,r,\cP)$ is the minimax rate of convergence on $N_{r,d}\left(\beta,L,\cP\right)$ when $1-\sum_{j=1}^d\frac{1}{\beta_j r_j}>0$ and that a proof of the corresponding lower bound can be obtained by a minor modification of that in Lepski and Willer \cite{lepski:lowerbound-deconvolution} to take into account the possible independence structure of the underlying density.

\section{Proofs of main results}
\label{proofs} 

\subsection{Quantities and technical lemma}
\label{sec:notations}
For brevity, introduce first
$$
\cI_d^{\diamond}:=\{I\in\cP\diamond\cP',\;(\cP,\cP')\in\mP\times\mP\},\quad\overline{\cU_p}:=\sup_{n\in\bN^*}\sup_{I\in\cI_d^{\diamond}}\sup_{h_I\in\cH_{p,I}}\cU_p(h_I)<\infty,
$$
Note that the finiteness of $\overline{\cU_p}$ is due both to the definition of the sets of multibandwidths $\cH_{p,I}$ and to the bounds given in Proposition \ref{multipliertheorem}.

\smallskip

Next, define the constant $\gamma_{p}$ involved in the selection rule. For $I\in\cI_d^{\diamond}$ and $\mathbf{r}\geq 1$, put
\begin{eqnarray*}
\gamma_{p,I}(\mathbf{r}):=\left\{
\begin{array}{llll}
4+\sqrt{\frac{37e^{-1}p\mathbf{r}}{2-p}}, &\ p\in(1,2),
\\*[4mm]
\left(7C_I+3\mathbf{A}\left(2\pi\right)^{-\frac{|I|}{2}}\left\|\widehat{K_{I}}g_I\right\|_{\infty}\left\|\widehat{q_I}\right\|_{\frac{1}{2}}\right)\mathbf{r}, &\ p=2,
\\*[4mm]
\left(\frac{46c(p)[p\vee e]}{3e}\right)c_{p}^{\frac{1}{p}-\frac{1}{2}}\left[1\vee C_I\right]\left(1\vee\left\|q_I\right\|_{\infty}\right)^{\frac{3}{4}}\mathbf{r}, &\ p\in(2,+\infty),
\\*[4mm]
6C_I(\mathbf{K},q)\left(1\vee\left\|q_I\right\|_{\infty}\right)^{\frac{1}{2}}\left[93|I|\ln(|I|)+69\mathbf{r}\right], &\ p=+\infty,
\end{array}
\right.
\end{eqnarray*}
where $g_I(t_I):=\prod_{j\in I}\left(1+t_j^2\right)^{\frac{\lambda_j}{2}}$, $C_I:=\mathbf{A}\left\{(2\pi)^{-\frac{|I|}{2}}\left(\left\|\widehat{K_I}g_I\right\|_{2}\vee\left\|\widehat{K_I}g_I\right\|_{1}\right)\right\}$, $c_p$ is given in the definition of $\cH_{p,I}$, $c(p):=15p/\ln(p)$ and
$$
C_I(\mathbf{K},q):=\frac{\mathbf{A}}{(2\pi)^{\frac{|I|}{2}}}\left\{\left\|\widehat{K_I}g_I\right\|_{2}\vee\left\|\widehat{K_I}g_I\right\|_{1}\vee\left(\max_{j\in I}\left\|D_j^1\widehat{K_I}g_I\right\|_{1}\right)\vee\left\|\widehat{K_I}\varphi_I\right\|_{2}\vee\left\|\widehat{K_I}\varphi_I\right\|_{1}\right\},
$$
with $\varphi_I(t_I):=\sup_{j\in I}\left|t_j\right|g_I(t_I)$. Then, put $\mathbf{r}_k:=kp\textbf{1}_{\{p<\infty\}}+k\textbf{1}_{\{p=+\infty\}}$, $k\geq 1$, and
\begin{eqnarray*}
\gamma_{p}:=\left\{
\begin{array}{ll}
\sup_{\cP,\cP'\in\mP}\sup_{I\in\cP\diamond\cP'}\left\{\gamma_{p,I}(\mathbf{r}_4)\right\}, &\ \mP\neq\{\overline{\emptyset}\},
\\*[4mm]
\gamma_{p,\overline{\emptyset}}(\mathbf{r}_1), &\ \mP=\{\overline{\emptyset}\}.
\end{array}
\right.
\end{eqnarray*}

\smallskip

Finally, we need the following technical lemma in order to compute our risk bounds. Define 
\begin{eqnarray*}
&&\xi_{p}:=\sup_{I\in\cI_d^{\diamond}}\sup_{h_I\in\cH_{p,I}}\left[\left\|\xi_{h_I}\right\|_{p}-\gamma_{p}\cU_{p}(h_I)\right]_+,
\\*[2mm]
&&\overline{\mathbf{f}}_{p}:=\md^2\left\|\textbf{K}\right\|_{1}^d\left[\overline{G}_{p}\right]^{\md(\md-1)}\left(\max\left\{\overline{G}_{p},\left\|\textbf{K}\right\|_{1}^d\mathbf{f}_p\right\}\right)^{\md-1},\quad\mathbf{f}_p:=1\vee\left[\sup_{I\in\cI_d^{\diamond}}\left\|f_I\right\|_{p}\right].
\end{eqnarray*}

\begin{lemma}
\label{lem:empiricalupperbound1} Assume that $\mP\neq\{\overline{\emptyset}\}$. Set $\mathbf{r}\in\{\mathbf{r}_1,\mathbf{r}_2,\mathbf{r}_4\}$. Under Assumptions (N1)-(N2), if $p\in(1,+\infty]$ then, for all integer $n\geq 3$,
$$
\left(\bE_f\left|\xi_{p}\right|^{\mathbf{r}}\right)^{\frac{1}{\mathbf{r}}}\leq\textbf{c}_{p,1}(\mathbf{r})n^{-\frac{1}{2}},\quad
\left(\bE_f\left|\overline{\mathbf{f}}_{p}\right|^{\mathbf{r}}\right)^{\frac{1}{\mathbf{r}}}\leq \textbf{c}_{p,2}(\mathbf{r},\mathbf{f}_p),\quad\forall f\in\bF_p\left[\;\mP\;\right].
$$
\end{lemma}
\noindent The absolute  constants $\textbf{c}_{p,1}(\mathbf{r})>0$ and $\textbf{c}_{p,2}(\mathbf{r},\mathbf{f}_p)>0$ can be explicitly expressed and the maps $\mathbf{f}_p\mapsto\textbf{c}_{p,2}(\mathbf{r},\mathbf{f}_p)$ are bounded on any bounded interval of $\bR_+$; see the proof of the latter result, which is postponed to Appendix.

\subsection{Oracle inequalities : proof of Theorems \ref{theo:oracleinequality1} and \ref{theo:oracleinequality2}.}
\label{sec:prooforacle} 

\paragraph{1)} Set $p\in[1,+\infty]$ and $f\in\bF_p\left[\;\mP\;\right]$. Let $(h,\cP)\in\cH_p\big[\;\mP\;\big]$, $\cP\in\mP(f)$, be fixed. 

\smallskip

In view of the triangle inequality we have
\begin{eqnarray*}
\left\|\widetilde{f}-f\right\|_p&\leq& \left\|\widetilde{f}_{(\widetilde{h},\widetilde{\cP})}-\widetilde{f}_{(h,\cP),(\widetilde{h},\widetilde{\cP})}\right\|_p+\left\|\widetilde{f}_{(h,\cP),(\widetilde{h},\widetilde{\cP})}-\widetilde{f}_{(h,\cP)}\right\|_p+\left\|\widetilde{f}_{(h,\cP)}-f\right\|_p
\\*[2mm]
&\leq& \widetilde{\Delta}_{p}(h,\cP)+\Lambda_{p}\cU_{p}(\widetilde{h},\widetilde{\cP})+\widetilde{\Delta}_{p}(\widetilde{h},\widetilde{\cP})+\Lambda_{p}\cU_{p}(h,\cP)+\left\|\widetilde{f}_{(h,\cP)}-f\right\|_p.
\end{eqnarray*}
Here we have used the equality $\widetilde{f}_{(h,\cP),(\widetilde{h},\widetilde{\cP})}=\widetilde{f}_{(\widetilde{h},\widetilde{\cP}),(h,\cP)}$. 
By definition of $(\widetilde{h},\widetilde{\cP})$, we obtain
\begin{equation}
\label{eq:globalloss1}
\left\|\widetilde{f}-f\right\|_p\leq 2\left[\widetilde{\Delta}_{p}(h,\cP)+\Lambda_{p}\cU_{p}(h,\cP)\right]+\left\|\widetilde{f}_{(h,\cP)}-f\right\|_p.
\end{equation}

\paragraph{2)}
Suppose that $\cP=\left\{I_1,\ldots,I_m\right\}$, $m\in\left\{1,\ldots,d\right\}$. Since $\cP\in\mP(f)$, for any $x\in\bR^d$ 
\begin{eqnarray*}
&&\left|\widetilde{f}_{(h,\cP)}(x)-f(x)\right|=\left|\prod_{I\in\cP}\widetilde{f}_{h_I}(x_{I})-\prod_{I\in\cP}f_I(x_{I})\right|
\\*[2mm]
&&\qquad\qquad\qquad\qquad\;\leq\sum_{j=1}^m\left|\widetilde{f}_{h_{I_j}}(x_{I_j})-f_{I_j}(x_{I_j})\right|\left(\prod_{k=\overline{j+1,m}}\left|\widetilde{f}_{h_{I_k}}(x_{I_k})\right|\right)\left(\prod_{l=\overline{1,j-1}}\left|f_{I_l}(x_{I_l})\right|\right).
\end{eqnarray*}
Here we have used the trivial equality: for $m\in\bN^*$ and $a_j,b_j\in\bR,\;j=\overline{1,m},$
\begin{equation}
\label{eq:productinequality}
\prod_{j=1}^{m}a_j-\prod_{j=1}^{m}b_j=\sum_{j=1}^m(a_j-b_j)\left(\prod_{k=\overline{j+1,m}}a_k\right)\left(\prod_{l=\overline{1,j-1}}b_l\right),
\end{equation}
where the product over empty set is assumed to be equal to one.

In view of $\cP\in\mP$, the triangle inequality and the Fubini-Tonelli theorem (used for the case $p<\infty$) we establish
\begin{eqnarray*}
&&\left\|\widetilde{f}_{(h,\cP)}-f\right\|_p\leq\sum_{j=1}^m\left\|\widetilde{f}_{h_{I_j}}-f_{I_j}\right\|_{p}\left(\prod_{k=\overline{j+1,m}}\left\|\widetilde{f}_{h_{I_k}}\right\|_{p}\right)\left(\prod_{l=\overline{1,j-1}}\left\|f_{I_l}\right\|_{p}\right)
\\*[2mm]
&&\qquad\qquad\qquad\leq m\left(\max\left\{\overline{G}_{p},\mathbf{f}_p\right\}\right)^{m-1}\sup_{I\in\cP}\left\|\widetilde{f}_{h_{I}}-f_{I}\right\|_{p},
\end{eqnarray*}
since $\left\|\textbf{K}\right\|_{1}\geq\int\textbf{K}=1$.
Remind that $\md=\sup_{\cP\in\mP}|\cP|$ and $\overline{G}_{p}\geq 1$. It follows
\begin{equation}
\label{eq:globalloss2}
\left\|\widetilde{f}_{(h,\cP)}-f\right\|_p\leq \md\left(\max\left\{\overline{G}_{p},\mathbf{f}_p\right\}\right)^{\md-1}\sup_{I\in\cP}\left\|\widetilde{f}_{h_{I}}-f_{I}\right\|_{p}.
\end{equation}

\paragraph{3)}
For any $(\eta,\cP')\in\cH_p\big[\;\mP\;\big]$ and any $x\in\bR^d$
\begin{eqnarray*}
\label{eq:globalloss3}
\left|\widetilde{f}_{(h,\cP),(\eta,\cP')}(x)-\widetilde{f}_{(\eta,\cP')}(x)\right|
=\left|\prod_{I'\in\cP'}\prod_{I\in\cP:I\cap I'\neq\emptyset}K_{\eta_{I\cap I'_j}}\star\widetilde{f}_{h_{I\cap I'_j}}(x_{I\cap I'})-\prod_{I'\in\cP'}\widetilde{f}_{\eta_{I'}}(x_{I'})\right|.
\end{eqnarray*}

Therefore, by the same method as the one used in step 2, we establish
\begin{equation}
\label{eq:globalloss4}
\left\|\widetilde{f}_{(h,\cP),(\eta,\cP')}-\widetilde{f}_{(\eta,\cP')}\right\|_p\leq \md\left[\;\overline{G}_{p}\;\right]^{\md(\md-1)}\sup_{I'\in\cP'}\left\|\prod_{I\in\cP:I\cap I'\neq\emptyset}\widetilde{f}_{h_{I\cap I'},\eta_{I\cap I'}}-\widetilde{f}_{\eta_{I'}}\right\|_{p}.
\end{equation}
Here we have used Young's inequality and the inequalities $\left\|\textbf{K}\right\|_{1}\geq\int\textbf{K}=1$ and $\overline{G}_{p}\geq 1$.

\paragraph{4)}
In view of the Fubini theorem and Young's inequality, for any $I\in\cI_d^{\diamond}$ and any $\eta\in(0,1]^d$
\begin{equation}
\label{eq:globalloss5}
\left\|\bE_f\left\{\widetilde{f}_{\eta_I}(\cdot)\right\}\right\|_{p}=\left\|K_{\eta_I}\star f_I\right\|_{p}\leq \left\|K_I\right\|_{1}\left\|f_I\right\|_{p}\leq\left\|\textbf{K}\right\|_{1}^d\mathbf{f}_p.
\end{equation}

Then, by the same method as the one used in step 2 and (\ref{eq:globalloss5}), for any $(\eta,\cP')\in\cH_p\big[\;\mP\;\big]$ and any $I'\in\cP'$ we get
\begin{eqnarray}
\label{eq:globalloss6}
&&\left\|\prod_{I\in\cP:I\cap I'\neq\emptyset}\widetilde{f}_{h_{I\cap I'},\eta_{I\cap I'}}-\prod_{I\in\cP:I\cap I'\neq\emptyset}\bE_f\left\{\widetilde{f}_{\eta_{I\cap I'}}(\cdot)\right\}\right\|_{p}\nonumber
\\*[2mm]
&&\leq \md\left(\max\left\{\overline{G}_{p},\left\|\textbf{K}\right\|_{1}^d\mathbf{f}_p\right\}\right)^{\md-1}\sup_{I\in\cP:I\cap I'\neq\emptyset}\left\|K_{\eta_{I\cap I'}}\star\left(\widetilde{f}_{h_{I\cap I'}}-f_{I\cap I'}\right)\right\|_{p}\nonumber
\\*[2mm]
&&\leq \md\left\|\textbf{K}\right\|_{1}^d\left(\max\left\{\overline{G}_{p},\left\|\textbf{K}\right\|_{1}^d\mathbf{f}_p\right\}\right)^{\md-1}\sup_{I\in\cP:I\cap I'\neq\emptyset}\left\|\widetilde{f}_{h_{I\cap I'}}-f_{I\cap I'}\right\|_{p}.
\end{eqnarray}

\paragraph{5)}
For $\eta\in(0,1]^d$ and $I'\in\cI_d$, since $\cP\in\mP(f)$, we have for any $x\in\bR^d$
$$
\bE_f\left\{\widetilde{f}_{\eta_{I'}}(x_{I'})\right\}=\int K_{\eta_{I'}}\left(y_{I'}-x_{I'}\right)\prod_{I\in\cP:I\cap I'\neq\emptyset}f_{I\cap I'}(y_{I\cap I'})\rd y_{I'}
=\prod_{I\in\cP:I\cap I'\neq\emptyset}\bE_f\left\{\widetilde{f}_{\eta_{I\cap I'}}(x_{I\cap I'})\right\}.
$$
Here we have used the product structure of the kernel $K$ and the Fubini theorem.

Thus, in view of the triangle inequality, (\ref{eq:globalloss4}), (\ref{eq:globalloss6}) and the trivial inequality $[\sup_ix_i-\sup_iy_i]_+\leq\sup_i[x_i-y_i]_+$, for any $(\eta,\cP')\in\cH_p\big[\;\mP\;\big]$, we get

\smallskip

\noindent $\left[\left\|\widetilde{f}_{(h,\cP),(\eta,\cP')}-\widetilde{f}_{(\eta,\cP')}\right\|_p-\Lambda_{p}\cU_{p}(\eta,\cP')\right]_+$
$$
\leq\md\left[\overline{G}_{p}\right]^{\md(\md-1)}\sup_{I'\in\cP'}\Bigg[\Big\|\prod_{I\in\cP:I\cap I'\neq\emptyset}\widetilde{f}_{h_{I\cap I'},\eta_{I\cap I'}}-\prod_{I\in\cP:I\cap I'\neq\emptyset}\bE_f\left\{\widetilde{f}_{\eta_{I\cap I'}}(\cdot)\right\}\Big\|_{p}+\left\|\xi_{\eta_{I'}}\right\|_{p}-\gamma_p\cU_p(\eta_{I'})\Bigg]_+;
$$

\smallskip

$$
\left[\left\|\widetilde{f}_{(h,\cP),(\eta,\cP')}-\widetilde{f}_{(\eta,\cP')}\right\|_p-\Lambda_{p}\cU_{p}(\eta,\cP')\right]_+\leq \overline{\mathbf{f}}_{p}\sup_{\cP'\in\mP}\sup_{I\in\cP\diamond\cP'}\left\|\widetilde{f}_{h_{I}}-f_{I}\right\|_{p}
+\overline{\mathbf{f}}_{p}\xi_{p},
$$
since $\overline{\mathbf{f}}_{p}\geq\md\left[\overline{G}_{p}\right]^{\md(\md-1)}\geq 1$. We deduce
\begin{equation}
\label{eq:globalloss7}
\widetilde{\Delta}_{p}(h,\cP)\leq\overline{\mathbf{f}}_{p}\sup_{\cP'\in\mP}\sup_{I\in\cP\diamond\cP'}\left\|\widetilde{f}_{h_{I}}-f_{I}\right\|_{p}
+\overline{\mathbf{f}}_{p}\xi_{p}.
\end{equation}
Finally, it follows from (\ref{eq:globalloss1}), (\ref{eq:globalloss2}) and (\ref{eq:globalloss7})
\begin{equation}
\label{eq:globalloss8}
\left\|\widetilde{f}-f\right\|_p\leq 3\overline{\mathbf{f}}_{p}\left\{\sup_{\cP'\in\mP}\sup_{I\in\cP\diamond\cP'}\left\|\widetilde{f}_{h_{I}}-f_{I}\right\|_{p}
+\gamma_p\cU_{p}(h,\cP)+\xi_{p}\right\}.
\end{equation}

\paragraph{6)} Consider the random event $B_p:=\Big\{\overline{G}_p\geq \cC_p\Big\}$, $\cC_p(\mathbf{f}_p)=\big(1+\gamma_p\overline{\cU_p}+\left\|\mathbf{K}\right\|_1^d\mathbf{f}_p\big)\left\|\mathbf{K}\right\|_1^d+1$.\\
\noindent Put also
$$
\cR_{p}^{(\mathbf{r})}\left[(h,\cP),f\right]:=\left(\bE_f\sup_{\cP'\in\mP}\sup_{I\in\cP\diamond\cP'}\left\|\widetilde{f}_{h_{I}}-f_{I}\right\|_{p}^{\mathbf{r}}\right)^{\frac{1}{\mathbf{r}}},\quad \mathbf{r}\geq 1.
$$

\smallskip

In view of (\ref{eq:globalloss5}), Lemma \ref{lem:empiricalupperbound1}, Markov's inequality, (\ref{eq:globalloss8}), and the Cauchy-Schwarz inequality we get $B_p\subseteq\big\{\xi_p\geq 1\big\}$, $\big[\bP_f\left(B_p\right)\big]^{\frac{1}{\mathbf{r}_4}}\leq\textbf{c}_{p,1}(\mathbf{r}_4)n^{-1/2}$ and
\begin{eqnarray*}
&&\left(\bE_f\left\|\widetilde{f}-f\right\|_p^{\mathbf{r}_1}1_{B_p^c}\right)^{\frac{1}{\mathbf{r}_1}}\leq 3\md^2\left\|\textbf{K}\right\|_{1}^d[\cC_p(\mathbf{f}_p)]^{\md^2-1}\Bigg(\cR_p^{(\mathbf{r}_1)}\left[(h,\cP),f\right]+\gamma_p\cU_p(h,\cP)+\frac{\textbf{c}_{p,1}(\mathbf{r}_1)}{\sqrt{n}}\Bigg),
\\*[2mm]
&&\left(\bE_f\left\|\widetilde{f}-f\right\|_p^{\mathbf{r}_1}1_{B_p}\right)^{\frac{1}{\mathbf{r}_1}}\leq 3\textbf{c}_{p,1}(\mathbf{r}_4)\textbf{c}_{p,2}(\mathbf{r}_4,\mathbf{f}_p)\left(\cR_p^{(\mathbf{r}_2)}\left[(h,\cP),f\right]+\gamma_p\overline{\cU_p}+\textbf{c}_{p,1}(\mathbf{r}_2)\right)n^{-1/2},
\\*[2mm]
&&\cR_p^{(\mathbf{r}_2)}\left[(h,\cP),f\right]\leq \textbf{c}_{p,1}(\mathbf{r}_2)+\gamma_p\overline{\cU_p}+\left\|\textbf{K}\right\|_{1}^d\mathbf{f}_p+\mathbf{f}_p.
\end{eqnarray*}
Thus, we come to the assertion of Theorem \ref{theo:oracleinequality1} with $\mathbf{C}_{p,1}(\mathbf{f}_p):=3\md^2\left\|\textbf{K}\right\|_{1}^d[\cC_p(\mathbf{f}_p)]^{\md^2-1}$ and 
\begin{eqnarray*}
&&\mathbf{C}_{p,2}(\mathbf{f}_p):=3\textbf{c}_{p,1}(\mathbf{r}_4)\textbf{c}_{p,2}(\mathbf{r}_4,\mathbf{f}_p)\left(2\gamma_p\overline{\cU_p}+(1+\left\|\textbf{K}\right\|_{1}^d)\mathbf{f}_p+2\textbf{c}_{p,1}{(\mathbf{r}_2)}\right)
\\*[2mm]
&&\qquad\qquad\qquad\qquad\qquad\qquad\qquad\qquad\qquad\qquad +3\textbf{c}_{p,1}(\mathbf{r}_1)\md^2\left\|\textbf{K}\right\|_{1}^d[\cC_p(\mathbf{f}_p)]^{\md^2-1},
\end{eqnarray*}
since $\cR_p^{(\mathbf{r}_1)}\left[(h,\cP),f\right]=\cR_p\left[(h,\cP),f\right]$. The constants $\textbf{c}_{p,1}(\mathbf{r}_k)$ and $\textbf{c}_{p,2}(\mathbf{r}_k,\mathbf{f}_p)$, $k=1,2,4$, are given in the proof of Lemma \ref{lem:empiricalupperbound1}.

\paragraph{7)} \textsl{\textbf{Particular case: $\mP=\big\{\overline{\emptyset}\big\}$ (no independence structure)}} 

\smallskip

Set $f\in\bF$ and let $h\in\cH_{p,\overline{\emptyset}}$ be arbitrary fixed. By scrutinizing the steps 1)-5) we easily see that
\begin{equation*}
\label{eq:globalloss9}
\left\|\widetilde{f}-f\right\|_p\leq (1+2\left\|\textbf{K}\right\|_{1}^d)\left\|\widetilde{f}_{h}-f\right\|_{p}
+2\gamma_{p,\overline{\emptyset}}(\mathbf{r}_1)\cU_{p}(h)+2\left[\left\|\xi_{h}\right\|_{p}-\gamma_{p,\overline{\emptyset}}(\mathbf{r}_1)\cU_{p}(h)\right]_+.
\end{equation*}
Thus, we get from Proposition \ref{prop:empiricalupperbound1}
\begin{equation*}
\label{eq:riskbound}
\left(\bE_f\left\|\widetilde{f}-f\right\|_p^{\mathbf{r}_1}\right)^{\frac{1}{\mathbf{r}_1}}\leq (1+2\left\|\textbf{K}\right\|_{1}^d)\left(\bE_f\left\|\widetilde{f}_h-f\right\|_p^{\mathbf{r}_1}\right)^{\frac{1}{\mathbf{r}_1}}
+2\gamma_{p,\overline{\emptyset}}(\mathbf{r}_1)\cU_{p}(h)+2c_p(\mathbf{r}_1) n^{-1/2},
\end{equation*}
where the constants $\gamma_{p,\overline{\emptyset}}(\mathbf{r}_1)$ and $c_p(\mathbf{r}_1)$ are given in the proof of Proposition \ref{prop:empiricalupperbound1}.\epr

\subsection{Adaptive minimax upper bounds: Proof of Theorems \ref{theo:Lpadaptiveupperbound}-\ref{theo:supadaptiveupperbound2}}

\paragraph{1)} \textsl{\textbf{Case $p\in(1,+\infty)$:}} 
let $\left(\beta,L,\cP\right)\in\left(0,l\right]^d\times(0,\infty)^d\times\mP$ and $f\in N_{p,d}\left(\beta,L,\cP\right)\subset\bF_p\left[\;\mP\;\right]$ be arbitrary fixed.

\smallskip 

In view of the triangle inequality, $\forall h\in (0,1]^d$,
\begin{equation}
\label{eq:upperbound11}
\sup_{\cP'\in\mP}\sup_{J\in\cP\diamond\cP'}\left\|\widetilde{f}_{h_J}-f_J\right\|_{p}\leq\sup_{\cP'\in\mP}\sup_{J\in\cP\diamond\cP'}\left\|\bE_f\{\widetilde{f}_{h_J}(\cdot)\}-f_J\right\|_{p}+\sup_{\cP'\in\mP}\sup_{J\in\cP\diamond\cP'}\left\|\xi_{h_J}\right\|_{p}
\end{equation}
where $\bE_f\{\widetilde{f}_{h_J}(x_J)\}=K_{h_J}\star f_J(x_J)$ and, remind, $\xi_{h_J}(x_J):=\widetilde{f}_{h_J}(x_J)-\bE_f\{\widetilde{f}_{h_J}(x_J)\}$.

\smallskip 

Note first that, by applying Proposition 3 in Kerkyacharian, Lepski and Picard \cite{kerk2}, it is easily established that, for any $h\in (0,1]^d$, any  $\cP'\in\mP$ and any $J\in\cP\diamond\cP'$,
\begin{equation}
\label{eq:upperbound13}
\left\|K_{h_J}\star f_J-f_J\right\|_{p}\leq \sum_{j\in J}c_J(\textbf{K},\left|J\right|,p,l,L_J)h_j^{\beta_j}\leq \textbf{c}\sum_{j\in J}h_j^{\beta_j}\leq\textbf{c}\sup_{I\in\cP}\sum_{j\in I}h_j^{\beta_j},\quad\textbf{c}>0.
\end{equation}

Next, if $(h,\cP)\in\cH_p[\;\mP\;]$, we easily get from Propositions \ref{multipliertheorem}-\ref{prop:empiricalupperbound1}
\begin{equation}
\label{eq:upperbound14}
\left(\bE_f\sup_{\cP'\in\mP}\sup_{J\in\cP\diamond\cP'}\left\|\xi_{h_J}\right\|_{p}^{p}\right)^{\frac{1}{p}}\leq \cO\left(\sup_{I\in\cP}\frac{1}{\sqrt{n\prod_{j\in I}h_j^{b_p^{-1}\lambda_j+1}}}\right).
\end{equation}

Consider now, for all $I\in\cP$, the system
$$
h_j^{\beta_j}=h_k^{\beta_k}=\frac{1}{\sqrt{n\prod_{j\in I}h_j^{b_p^{-1}\lambda_j+1}}},\quad j,k\in I.
$$
The solution is given by
\begin{equation}
\label{eq:adaptiveminimaxbandwidth11}
h_j=n^{-\frac{\tau_I}{2\tau_I+1}\frac{1}{\beta_j}},\quad j\in I, \quad I\in\cP,
\end{equation}
where $\tau_I$ is given in (\ref{eq:Lpminimaxrate}).

\smallskip

Note that, for all $I\in\cP$, $h_I\in[h_{min}^{(p)},h_{max}^{(p)}]^{|I|}$ and $n\prod_{j\in I}h_j^{b_p^{-1}\lambda_j+1}\geq 1$ for $n$ large enough. Denote by $\overline{h}_I$ the projection of $h_I$ on the dyadic grid $\cH_{p,I}$. It is easily checked that $(\overline{h},\cP)\in\cH_p\big[\;\mP\;\big]$ for $n$ large enough. Thus, it follows from Theorem \ref{theo:oracleinequality1}, (\ref{eq:upperbound11}), (\ref{eq:upperbound13}) and (\ref{eq:upperbound14}) that
 
\begin{equation}
\label{eq:adaptiveupperbound2}
\cR_p\left[\widehat{f},f\right]\leq \textbf{C}\left[\sup_{I\in\cP}\sum_{j\in I}\overline{h}_j^{\beta_j}+\sup_{I\in\cP}\frac{1}{\sqrt{n\prod_{j\in I}\overline{h}_j^{b_p^{-1}\lambda_j+1}}}\right]+\alpha_{p,2}n^{-1/2},\quad\textbf{C}>0,
\end{equation}
for $n$ large enough. Finally, it is easily seen that we get the statement of Theorem \ref{theo:Lpadaptiveupperbound} from (\ref{eq:adaptiveminimaxbandwidth11}) and (\ref{eq:adaptiveupperbound2}). Similarly, the statement of Theorem \ref{theo:Lpadaptiveupperbound2} is obtained by applying Theorem \ref{theo:oracleinequality2}.

\paragraph{2)} \textsl{\textbf{Case $p=+\infty$:}} 
let $\left(\beta,L,r,\cP\right)\in\left(0,l\right]^d\times(0,\infty)^d\times[1,+\infty]^d\times\mP$ such that $1-\sum_{j=1}^d\frac{1}{\beta_j r_j}>0$ and $f\in N_{r,d}\left(\beta,L,\cP\right)$ be arbitrary fixed. It follows from the definition of the latter functional class and the embedding theorem for anisotropic Nikolskii classes, see, e.g., Theorem 6.9 in Nikolskii \cite{nikolskii}, that $N_{r,d}\left(\beta,L,\cP\right)\subset\bF_{\infty}\left[\;\mP\;\right]$, since $1-\sum_{j\in I}\frac{1}{\beta_j r_j}>0$, $\forall I\in\cI_d$.

\smallskip 

Note first that, in view of the arguments given in the proof of Theorem 3 in Lepski \cite{lepski:supnormlossdensityestimation}, it follows from Lemme 4 in the latter paper that, for any $h\in (0,1]^d$, any  $\cP'\in\mP$ and any $J\in\cP\diamond\cP'$,
\begin{eqnarray}
\label{eq:supperbound5}
&&\left\|K_{h_J}\star f_J-f_J\right\|_{\infty}\leq\textbf{c}\sup_{I\in\cP}\sum_{j\in I}h_j^{\beta_j(I)},\quad\textbf{c}:=c(\textbf{K},d,l,L)>0,
\\*[2mm]
&&\beta_j(I):=\sigma(I)\beta_i\sigma_j^{-1}(I),\quad\sigma(I):=1-\sum_{k\in I}\left(\beta_kp_k\right)^{-1},\quad\sigma_j(I):=1-\sum_{k\in I}\left(p_k^{-1}-p_j^{-1}\right)\beta_k^{-1}.\nonumber
\end{eqnarray}

Next, if $(h,\cP)\in\cH_{\infty}[\;\mP\;]$, we easily get from Propositions \ref{prop:empiricalupperbound1}
\begin{equation}
\label{eq:supperbound6}
\bE_f\sup_{\cP'\in\mP}\sup_{J\in\cP\diamond\cP'}\left\|\xi_{h_I}\right\|_{\infty}\leq \cO\left(\sup_{I\in\cP}\sqrt{\frac{\ln(n)}{n\prod_{j\in I}h_j^{2\lambda_j+1}}}\right).
\end{equation}

Consider now, for all $I\in\cP$, the system
$$
h_j^{\beta_j(I)}=h_k^{\beta_k(I)}=\sqrt{\frac{\ln(n)}{n\prod_{j\in I}h_j^{2\lambda_j+1}}},\quad j,k\in I.
$$
The solution is given by
\begin{equation}
\label{eq:adaptiveminimaxbandwidth1}
h_j=\left(\frac{n}{\ln(n)}\right)^{-\frac{\Upsilon_I}{2\Upsilon_I+1}\frac{1}{\beta_j(I)}},\quad j\in I, \quad I\in\cP,
\end{equation}
where $\Upsilon_I$ is given in (\ref{eq:suprate}).

\smallskip

Note that, for all $I\in\cP$, $n\prod_{j\in I}h_j^{2\lambda_j+1}\geq \ln(n)$ for $n$ large enough. Thus, as previously, we get the statement of Theorem \ref{theo:supadaptiveupperbound} from Theorem \ref{theo:oracleinequality1}, (\ref{eq:supperbound5}), (\ref{eq:supperbound6}) and (\ref{eq:adaptiveminimaxbandwidth1}). Similarly, the statement of Theorem \ref{theo:supadaptiveupperbound2} is obtained by applying Theorem \ref{theo:oracleinequality2}.
\epr

\section{Appendix}

\subsection{Proof of Proposition \ref{multipliertheorem}} Let $h\in(0,1]^d$ and $I\in\cI_d^{\diamond}$ be arbitrary fixed. Note that
\begin{equation}
\label{eq:auxiliaryfunction}
L_{(h_I)}(x_I):=\frac{1}{(2\pi)^{\left|I\right|}}\int_{\bR^{\left|I\right|}}e^{-i\left\langle t_I,x_I\right\rangle}\frac{\widehat{K_{I}}(h_It_I)g_{I}(h_It_I)}{g_{I}(h_It_I)\widehat{q_I}(t_I)}\rd t_I,\quad g_I(t_I):=\prod_{j\in I}\left(1+t_j^2\right)^{\frac{\lambda_j}{2}},
\end{equation}
where $h_I t_I$ denotes the coordinate-wise product of the vectors $h_I$ and $t_I$.

\paragraph{1)\textit{Proof of assertion (i)}} Set $p\in(1,2)$. Here, we apply the Marcinkiewicz Multiplier Theorem on $\bR^{|I|}$, given in Grafakos \cite{Grafakos:multiplier} p. 363, with
$$
m(t_I)=g_I^{-1}(h_It_I)\widehat{q_I}^{-1}(t_I).
$$

In view of Assumption (N2) on $q$, $m$ is a bounded function defined away from the coordinates axes on $\bR^{|I|}$ and is $\cC^{|I|}$ on this region. Moreover,
\begin{equation}
\label{eq:multiplier}
\sup_{t_I\in\bR^{|I|}}\left|m(t_I)\right|\leq\mathbf{A}\sup_{u_I\in\bR^{|I|}}\left[\prod_{j\in I}\left(1+u_j^2\right)^{-\frac{\lambda_j}{2}}\prod_{j\in I}\left(1+[u_j/h_j]^2\right)^{\frac{\lambda_j}{2}}\right]\leq\mathbf{A}\prod_{j\in I}h_j^{-\lambda_j}.
\end{equation}

Set $\alpha_I=(\alpha_j)_{j\in I}\in\bN^{|I|}$ satisfying $|\alpha_I|:=\sum_{j\in I}\alpha_j\leq|I|$. In view of Leibniz's rule, one has
\begin{eqnarray*}
\left[D^{\alpha_I}m\right](t_I)=\sum_{\gamma_I\leq\alpha_I}\left(\begin{array}{cc}
\alpha_I\\
\gamma_I
\end{array}\right)
\left\{\prod_{j\in I}h_j^{\gamma_j}\right\}\left[D^{\gamma_I}g_I^{-1}\right](h_I t_I)\left[D^{\alpha_I-\gamma_I}\widehat{q_I}^{-1}\right](t_I),\quad\forall t\in\bR^{d}.
\end{eqnarray*}
Here, $\gamma_I\leq\alpha_I$ means $\gamma_j\leq\alpha_j$, $\forall j\in I$.

\smallskip

Let $t_I$ be chosen such that $t_j\neq 0$ if $\alpha_j\neq 0$. In this case, for any multi-index $\gamma_I\leq\alpha_I$,
$$
\left\{\prod_{j\in I}h_j^{\gamma_j}\right\}\left[D^{\gamma_I}g_I^{-1}\right](h_I t_I)\left[D^{\alpha_I-\gamma_I}\widehat{q_I}^{-1}\right](t_I)=
$$
$$
\left\{\prod_{j\in I}\left(t_jh_j\right)^{\gamma_j}\right\}\left[D^{\gamma_I}g_I^{-1}\right](h_I t_I)\left[D^{\alpha_I-\gamma_I}\widehat{q_I}^{-1}\right](t_I)\left\{\prod_{j\in I}t_j^{\alpha_j-\gamma_j}\right\}\left(\prod_{j\in I}t_j^{-\alpha_j}\right).
$$
Here, we assume that $0^0$ is equal to one.

\smallskip

Since $q$ satisfies Assumption (N2), we obtain similarly as in (\ref{eq:multiplier})
\begin{eqnarray}
\label{eq:multiplier2}
&&\left|\;\left[D^{\alpha_I}m\right](t_I)\right|\leq C(|I|,q_I)\mathbf{A}\left\{\prod_{j\in I}h_j^{-\lambda_j}\right\}\left(\prod_{j\in I}\left|t_j\right|^{-\alpha_j}\right),
\\*[2mm]
&& C(|I|,q_I):=\max_{|\alpha_I|\leq|I|}\left\{\sum_{\gamma_I\leq\alpha_I}\left(\begin{array}{cc}
\alpha_I\\
\gamma_I
\end{array}\right)\sup_{u_I\in\bR^{|I|}}\left|\left\{\prod_{j\in I}u_j^{\gamma_j}\right\}\left[D^{\gamma_I}g_I^{-1}\right](u_I)g_I(u_I)\right|\right\}<\infty.
\end{eqnarray}

Put $\widehat{S_{I}}(t_I):=\widehat{K_{I}}(t_I)g_{I}(t_I)$, $t\in\bR^{d}$. Since $\mathbf{K}\in\bS(\bR)$, $\widehat{S_{I}}\in\bS(\bR^{|I|})$ is the Fourier transform of a function $S_{I}\in\bS(\bR^{|I|})\subset\bL_p(\bR^{|I|})$. As
$$
L_{(h_I)}(x_I):=\frac{1}{(2\pi)^{\left|I\right|}}\int_{\bR^{\left|I\right|}}e^{-i\left\langle t_I,x_I\right\rangle}m(t_I)\widehat{S_{I}}(h_It_I)\rd t_I,
$$
it follows from Corollary 5.2.5. in Grafakos \cite{Grafakos:multiplier}, (\ref{eq:multiplier}) and (\ref{eq:multiplier2})
$$
\left\|L_{(h_I)}\right\|_{p}\leq 2\mathbf{A}C_{|I|}C(|I|,q_I)\max\left(p,(p-1)^{-1}\right)^{6|I|}\left\|S_{I}\right\|_{p}\left(V_{h_I}\right)^{-(1-1/p)}\prod_{j\in I}h_j^{-\lambda_j},
$$
where $C_{|I|}<\infty$ is a dimensional constant which is not explicitly done in the aforementioned result. Thus, assertion (i) of Proposition \ref{multipliertheorem} is proved with
$$
C_{p,I}:=2\mathbf{A}\left\{C_{|I|}C(|I|,q_I)\max\left(p,(p-1)^{-1}\right)^{6|I|}\left\|S_{I}\right\|_{p}\right\}.
$$

\smallskip

\paragraph{2)\textit{Proof of assertion (ii)}} Note first that
$$
\left\|L_{(h_I)}\right\|_{2}=(2\pi)^{-\frac{|I|}{2}}\left\|\widehat{K_{h_I}}/\widehat{q_I}\right\|_{2},\quad\left\|L_{(h_I)}\right\|_{\infty}\leq(2\pi)^{-|I|}\left\|L_{(h_I)}\right\|_{1}.
$$
In view of Assumption (N2) on the errors,
\begin{eqnarray*}
\left\|\widehat{K_{h_I}}/\widehat{q_I}\right\|_{2}^2 &\leq& \mathbf{A}^2\int_{\bR^{|I|}}\left|\widehat{K_{I}}(h_It_I)\right|^2\prod_{j\in I}(1+t_j^2)^{\lambda_j}\rd t_I
\\*[2mm]
&\leq& \mathbf{A}^2\left(\int_{\bR^{|I|}}\left|\widehat{K_{I}}(u_I)\right|^2\prod_{j\in I}(1+u_j^2)^{\lambda_j}\rd u_I\right)V_{h_I}^{-1}\prod_{j\in I}h_j^{-2\lambda_j};
\\*[2mm]
\left\|\widehat{K_{h_I}}/\widehat{q_I}\right\|_{1} &\leq& \mathbf{A}\left(\int_{\bR^{|I|}}\left|\widehat{K_{I}}(u_I)\right|\prod_{j\in I}(1+u_j^2)^{\lambda_j/2}\rd u_I\right)V_{h_I}^{-1}\prod_{j\in I}h_j^{-\lambda_j}.
\end{eqnarray*}

Thus, assertion (ii) of Proposition \ref{multipliertheorem} is proved with
$$
C_I:=\mathbf{A}\left\{(2\pi)^{-\frac{|I|}{2}}\left(\left\|\widehat{K_I}g_I\right\|_{2}\vee\left\|\widehat{K_I}g_I\right\|_{1}\right)\right\},
$$
where $g_I$ is given in (\ref{eq:auxiliaryfunction}).
\epr

\subsection{Proof of Proposition \ref{prop:empiricalupperbound1}: case $p<\infty$.} Let $I\in\cI_d^{\diamond}$ be arbitrary fixed. We get the statement of Proposition \ref{prop:empiricalupperbound1} by applying Theorem 1 and Corollaries 2 and 3 in Goldenshluger and Lepski \cite{G-L:uniformbounds} with $s=p$, $\cX=\cT=\bR^{|I|}$, $\nu=\tau$ is the Lebesgue measure on $\bR^{|I|}$, $w(\cdot,\cdot)=n^{-1}L_{(h_I)}(\cdot-\cdot)$ and $M_s(w)=\left\|n^{-1}L_{(h_I)}\right\|_{p}<\infty$. Here, the i.i.d. random vectors are the $Y_{k,I}$'s and their common density is $f_I\star q_I$. By using the continuity property of $L_{(h_I)}(\cdot)$, it is easily proved that Assumption (A1) in the aforementioned paper is fulfilled. 

\smallskip

\noindent\textbf{1)} \textsl{Case $p\in(1,2)$}. Let $\mathbf{r}\geq 1$ and $h_I\in\cH_{p,I}$ be arbitrary fixed. 

\smallskip

By application of Corollary 2 in Goldenshluger and Lepski \cite{G-L:uniformbounds}, one has
\begin{equation}
\label{eq:uniformbound1}
\bP\left\{\left\|\xi_{h_I}\right\|_{p}\geq U_{p}(h_I)+z\right\}\leq\exp\left\{-\frac{z^2}{A_p^2(h_I)}\right\},\quad \forall z>0,\;\forall n\geq 1,
\end{equation}
where $U_{p}(h_I)=4n^{\frac{1}{p}-1}\left\|L_{(h_I)}\right\|_{p}$ and $A_p^2(h_I)=37n^{-1}\left\|L_{(h_I)}\right\|_{p}^2$. 

\smallskip

By integration of (\ref{eq:uniformbound1}) we easily get, for all integer $n\geq 3$,
\begin{eqnarray*}
\label{eq:uniformbound2}
\bE\left[\left\|\xi_{h_I}\right\|_{p}-U_{p}(h_I)-A_p(h_I)\sqrt{r\ln(n)}\right]_+^{\mathbf{r}}
&\leq& \Gamma(\mathbf{r}+1)\left[A_p(h_I)\right]^{\mathbf{r}} e^{-\mathbf{r}\ln(n)}
\\*[2mm]
&\leq& \Gamma(\mathbf{r}+1)7^{\mathbf{r}}\sup_{h_I\in\cH_{p,I}}\left[n^{-\frac{1}{2}}\left\|L_{(h_I)}\right\|_{p}\right]^{\mathbf{r}} n^{-\mathbf{r}},
\end{eqnarray*}
where $\Gamma(\cdot)$ is the well known Gamma function.

\smallskip

Note that, for all integer $n\geq 3$, 
\begin{eqnarray*}
U_{p}(h_I)+A_p(h_I)\sqrt{\mathbf{r}\ln(n)}\leq \left\{4+\sqrt{\frac{37e^{-1}p\mathbf{r}}{2-p}}\;\right\}n^{\frac{1}{p}-1}\left\|L_{(h_I)}\right\|_{p}=:\gamma_{p,I}(\mathbf{r})\cU_p(h_I).
\end{eqnarray*}

\smallskip

Since card$(\cH_{p,I})\leq\Big[\Big(1\vee\frac{p}{|I|}\Big)\log_2(n)\Big]^{|I|}$, we obtain, for all integer $n\geq 3$,
\begin{eqnarray*}
&&\left\{\bE\sup_{h_I\in\cH_{p,I}}\left[\left\|\xi_{h_I}\right\|_{p}- \gamma_{p,I}(\mathbf{r})\cU_{p}(h_I)\right]_+^{\mathbf{r}}\right\}^{\frac{1}{\mathbf{r}}}\leq c_p(\mathbf{r}) n^{-\frac{1}{2}},
\\*[2mm]
&& c_p(\mathbf{r}):=7\left[\Gamma(\mathbf{r}+1)\right]^{\frac{1}{\mathbf{r}}}\sup_{n\in\bN^*}\sup_{I\in\cI_d^{\diamond}}\sup_{h_I\in\cH_{p,I}} \left\{n^{-1}[2\log_2(n)]^{\frac{|I|}{\mathbf{r}}}\left\|L_{(h_I)}\right\|_{p}\right\},
\end{eqnarray*}
which is finite in view of Proposition \ref{multipliertheorem} and the definition of the set $\cH_{p,I}$.

\smallskip

\noindent\textbf{2)} \textsl{Case $p=2$}. Let $\mathbf{r}\geq 1$ and $h_I\in\cH_{p,I}$ be arbitrary fixed. Here, we apply Theorem 1 in Goldenshluger and Lepski \cite{G-L:uniformbounds} but we compute differently the upper bound on the "\textit{dual}" variance $\sigma^2$ by using the arguments given in the proof of Proposition 7 in Comte and Lacour \cite{comte:deconvolution}. Indeed, we obtain
$$
\sigma^2\leq n^{-2}\left(2\pi\right)^{-|I|}\left\|\frac{\widehat{K_{h_I}}}{\widehat{q_I}}\right\|_{\infty}^2\int_{\bR^{|I|}}\left|\widehat{f_I}(t_I)\widehat{q_I}(t_I)\right|\rd t_I\leq n^{-2}\left(2\pi\right)^{-|I|}\left\|\widehat{q_I}\right\|_{1}\left\|\frac{\widehat{K_{h_I}}}{\widehat{q_I}}\right\|_{\infty}^2,
$$
since $\big\|\widehat{f_I}\big\|_{\infty}\leq\left\|f_I\right\|_{1}=1$. 

\smallskip

Taking into account the latter inequality, the result of Theorem 1 in Goldenshluger and Lepski \cite{G-L:uniformbounds} should be
\begin{eqnarray}
\label{eq:uniformbound3}
&&\bP\left\{\left\|\xi_{h_I}\right\|_{p}\geq U_{p}(h_I)+z\right\}\leq\exp\left\{-\frac{z^2}{A_p^2(h_I)+B_p(h_I)z}\right\},\quad \forall z>0,\;\forall n\geq 1,
\\*[2mm]
&& U_{p}(h_I)=n^{-\frac{1}{2}}\left\|L_{(h_I)}\right\|_{2},\nonumber
\\*[2mm]
&& A_p^2(h_I)=\frac{6}{\left(2\pi\right)^{|I|}}\left\|\widehat{q_I}\right\|_{1}n^{-1}\left\|\widehat{K_{h_I}}/\widehat{q_I}\right\|_{\infty}^2+24n^{-\frac{3}{2}}\left\|L_{(h_I)}\right\|_{2}^2,\;\; B_p(h_I)=\frac{4}{3}n^{-1}\left\|L_{(h_I)}\right\|_{2}.\nonumber
\end{eqnarray}

By integration of (\ref{eq:uniformbound3}) we get, for all integer $n\geq 3$,
\begin{eqnarray*}
&&\bE\left[\left\|\xi_{h_I}\right\|_{p}- U_{p}(h_I)-A_p(h_I)\sqrt{ \mathbf{r}\ln(n)}-B_p(h_I) \mathbf{r}\ln(n)\right]_+^{\mathbf{r}}
\\*[2mm]
&&\leq\Gamma(\mathbf{r}+1)\left\{A_p(h_I)+B_p(h_I)\right\}^{\mathbf{r}} e^{- \mathbf{r}\ln(n)}
\\*[2mm]
&&\leq\Gamma(\mathbf{r}+1)\left(6\vee\left\|\widehat{q_I}\right\|_{1}^{\frac{1}{2}}\right)^{\mathbf{r}}\sup_{h_I\in\cH_{p,I}} \left\{\left\|\widehat{K_{h_I}}/\widehat{q_I}\right\|_{\infty}+\left\|L_{(h_I)}\right\|_{2}\right\}^{\mathbf{r}} n^{-\frac{\mathbf{r}}{2}-\mathbf{r}}.
\end{eqnarray*}

Note that, in view of Assumption (N2) on the errors,
\begin{equation}
\label{eq:supnorm}
\left\|\widehat{K_{h_I}}/\widehat{q_I}\right\|_{\infty}\leq\mathbf{A}\left\|\widehat{K_{I}}g_I\right\|_{\infty}\prod_{j\in I}h_j^{-\lambda_j},
\end{equation}
where $g_I$ is given in (\ref{eq:auxiliaryfunction}). Thus, in view of Proposition \ref{multipliertheorem}, (\ref{eq:supnorm}) and the definition of $\cH_{p,I}$, for all integer $n\geq 3$, 
\begin{eqnarray*}
&& U_{p}(h_I)+A_p(h_I)\sqrt{ \mathbf{r}\ln(n)}+B_p(h_I) \mathbf{r}\ln(n)
\\*[2mm]
&&\leq \left\{C_I\left(1+\frac{8\mathbf{r}}{3e}+\sqrt{\frac{48\mathbf{r}}{e}}\right)+\mathbf{A}\left\|\widehat{K_{I}}g_I\right\|_{\infty}\sqrt{\frac{6\mathbf{r}\left\|\widehat{q_I}\right\|_{1}}{\left(2\pi\right)^{|I|}}}\right\}n^{-\frac{1}{2}}\prod_{j\in I}h_j^{-\lambda_j-\frac{1}{2}}
\\*[2mm]
&& =:\gamma_{p,I}(\mathbf{r})\cU_p(h_I).
\end{eqnarray*}

Finally, we obtain for all integer $n\geq 3$
\begin{eqnarray*}
&&\left\{\bE\sup_{h_I\in\cH_{p,I}}\left[\left\|\xi_{h_I}\right\|_{p}-\gamma_{p,I}(\mathbf{r}) \cU_{p}(h_I)\right]_+^{\mathbf{r}}\right\}^{^{\frac{1}{\mathbf{r}}}}\leq c_p(\mathbf{r}) n^{-\frac{1}{2}},\quad c_p(\mathbf{r}):=\left[\Gamma(\mathbf{r}+1)\right]^{\frac{1}{\mathbf{r}}}
\\*[2mm]
&&\qquad\qquad\times\sup_{n\in\bN^*}\sup_{I\in\cI_d^{\diamond}}\sup_{h_I\in\cH_{p,I}} \left[\left(6\vee\left\|\widehat{q_I}\right\|_{1}^{\frac{1}{2}}\right)\left\{\left\|\widehat{K_{h_I}}/\widehat{q_I}\right\|_{\infty}+\left\|L_{(h_I)}\right\|_{2}\right\}\big[2\log_2(n)\big]^{\frac{|I|}{\mathbf{r}}}n^{-1}\right],
\end{eqnarray*} 
which is finite in view of Proposition \ref{multipliertheorem}, (\ref{eq:supnorm}) and the definition of the set $\cH_{p,I}$.

\smallskip

\noindent\textbf{3)} \textsl{Case $p>2$}. Let $\mathbf{r}\geq 1$ and $h_I\in\cH_{p,I}$ be arbitrary fixed. 

\smallskip

By application of Corollary 3 in Goldenshluger and Lepski \cite{G-L:uniformbounds}, one has
\begin{eqnarray}
\label{eq:uniformbound4}
&&\bP\left\{\left\|\xi_{h_I}\right\|_{p}\geq U_{p}(h_I)+z\right\}\leq\exp\left\{-\frac{z^2}{A_p^2(h_I)+B_p(h_I)z}\right\},\quad \forall z>0,\;\forall n\geq 1,
\\*[2mm]
&& U_{p}(h_I)=3c(p)\left\|q_I\right\|_{\infty}^{\frac{1}{2}-\frac{1}{p}}\left\{n^{-\frac{1}{2}}\left\|L_{(h_I)}\right\|_{2}+n^{\frac{1}{p}-1}\left\|L_{(h_I)}\right\|_{p}\right\},\nonumber
\end{eqnarray}
\begin{eqnarray*} A_p^2(h_I)&=&16c(p)\left\|q_I\right\|_{\infty}^{\frac{3}{2}}\left\{n^{-1}\left\|L_{(h_I)}\right\|_{\frac{2p}{p+2}}^2+n^{-\frac{3}{2}}\left\|L_{(h_I)}\right\|_{2}\left\|L_{(h_I)}\right\|_{p}+n^{\frac{1}{p}-2}\left\|L_{(h_I)}\right\|_{p}^2\right\},\nonumber
\\*[2mm]
B_p(h_I)&=&\frac{4}{3}c(p)n^{-1}\left\|L_{(h_I)}\right\|_{p},\quad c(p)=\frac{15p}{\ln(p)}.\nonumber
\end{eqnarray*}
Here, we have used the following inequalities, which are consequences of Young's inequality.
\begin{eqnarray*}
&&\left\|f_I\star q_I\right\|_{\infty}\leq\left\|f_I\right\|_{1}\left\|q_I\right\|_{\infty}\leq \left\|q_I\right\|_{\infty},
\\*[2mm]
&&\left\|\sqrt{f_I\star q_I}\right\|_{p}\leq\left\|f_I\star q_I\right\|_{\infty}^{\frac{1}{2}-\frac{1}{p}}\left\|f_I\star q_I\right\|_{1}\leq \left\|q_I\right\|_{\infty}^{\frac{1}{2}-\frac{1}{p}}.
\end{eqnarray*}

By integration of (\ref{eq:uniformbound4}) we get, for all integer $n\geq 3$,
\begin{eqnarray*}
&&\bE\left[\left\|\xi_{h_I}\right\|_{p}- U_{p}(h_I)-A_p(h_I)\sqrt{\mathbf{r}\ln(n)}-B_p(h_I)\mathbf{r}\ln(n)\right]_+^{\mathbf{r}}
\\*[2mm]
&&\leq\Gamma(\mathbf{r}+1)\left\{A_p(h_I)+B_p(h_I)\right\}^{\mathbf{r}} e^{-\mathbf{r}\ln(n)}
\\*[2mm]
&&\leq\Gamma(\mathbf{r}+1)\left\{6c(p)\left(1\vee\left\|q_I\right\|_{\infty}^{\frac{3}{4}}\right)\right\}^{\mathbf{r}}
\\*[2mm]
&&\qquad\qquad\times\sup_{h_I\in\cH_{p,I}} \left\{\left\|L_{(h_I)}\right\|_{\frac{2p}{p+2}}+\sqrt{\left\|L_{(h_I)}\right\|_{2}\left\|L_{(h_I)}\right\|_{p}}+\left\|L_{(h_I)}\right\|_{p}\right\}^{\mathbf{r}} n^{-\frac{\mathbf{r}}{2}-\mathbf{r}}.
\end{eqnarray*}

\smallskip

In view of Proposition \ref{multipliertheorem}, we get
\begin{equation}
\label{eq:Lp}
\left\|L_{(h_I)}\right\|_{p}\leq\left\|L_{(h_I)}\right\|_{\infty}^{1-\frac{2}{p}}\left\|L_{(h_I)}\right\|_{2}^{\frac{2}{p}}\leq C_IV_{h_I}^{\frac{1}{p}-\frac{1}{2}}\prod_{j\in I}h_j^{-\lambda_j-\frac{1}{2}}.
\end{equation}

Thus, in view of Proposition \ref{multipliertheorem}, (\ref{eq:Lp}) and the definition of $\cH_{p,I}$, for all integer $n\geq 3$, 
\begin{eqnarray*}
&& U_{p}(h_I)+A_p(h_I)\sqrt{\mathbf{r}\ln(n)}+B_p(h_I)\mathbf{r}\ln(n)
\\*[2mm]
&&\leq c_{p}^{\frac{1}{p}-\frac{1}{2}}\left(1\vee C_I\right)\left\{6c(p)\left\|q_I\right\|_{\infty}^{\frac{1}{2}-\frac{1}{p}}+8\sqrt{\frac{\mathbf{r}c(p)[p\vee e]}{e}}\left\|q_I\right\|_{\infty}^{\frac{3}{4}}+\frac{4\mathbf{r}c(p)[p\vee e]}{3e}\right\}
\\*[2mm]
&&\qquad\times n^{-\frac{1}{2}}\left[\prod_{j\in I}h_j^{-\lambda_j-\frac{1}{2}}+\sqrt{\ln(n)}\left\|L_{(h_I)}\right\|_{\frac{2p}{p+2}}\right]
\\*[2mm]
&& :=\gamma_{p,I}(\mathbf{r})\cU_p(h_I).
\end{eqnarray*}

\smallskip

Finally, we obtain for all integer $n\geq 3$
\begin{eqnarray*}
&&\left\{\bE\sup_{h_I\in\cH_{p,I}}\left[\left\|\xi_{h_I}\right\|_{p}-\gamma_{p,I}(\mathbf{r}) \cU_{p}(h_I)\right]_+^{\mathbf{r}}\right\}^{^{\frac{1}{\mathbf{r}}}}\leq c_p(\mathbf{r}) n^{-\frac{1}{2}},\qquad c_p(\mathbf{r}):=6c(p)\left[\Gamma(\mathbf{r}+1)\right]^{\frac{1}{\mathbf{r}}}
\\*[2mm]
&& \times\sup_{n\in\bN^*}\sup_{I\in\cI_d^{\diamond}}\sup_{h_I\in\cH_{p,I}} \Big[\;\Big(1\vee\left\|q_I\right\|_{\infty}^{\frac{3}{4}}\Big)\Big\{\left\|L_{(h_I)}\right\|_{\frac{2p}{p+2}}+\sqrt{\left\|L_{(h_I)}\right\|_{2}\left\|L_{(h_I)}\right\|_{p}}
\\*[2mm]
&&\qquad\qquad\qquad\qquad\qquad\qquad\qquad\qquad\qquad\qquad+\left\|L_{(h_I)}\right\|_{p}\Big\}n^{-1}\Big[\Big(1\vee\frac{p}{|I|}\Big)\log_2(n)\Big]^{\frac{|I|}{\mathbf{r}}}\;\Big],
\end{eqnarray*} 
which is finite in view of Proposition \ref{multipliertheorem}, (\ref{eq:Lp}) and the definition of the set $\cH_{p,I}$.
\epr

\subsection{Proof of Proposition \ref{prop:empiricalupperbound1}: case $p=+\infty$.} Let $n\geq 3$, $I\in\cI_d^{\diamond}$ and $h_I\in[1/n,1]^{|I|}$ be arbitrary fixed. Assume that $n\prod_{j\in I}h_j^{2\lambda_j+1}\geq\ln(n)$. We divide this proof into several steps.

\paragraph{1) Preliminaries: } First, since $q$ satisfies Assumption (N2) and the $Y_{k,I}$'s are i.i.d. random vectors with density $f_I\star q_I$, we get from Proposition \ref{multipliertheorem}
\begin{eqnarray}
\label{eq:uniformbound11}
&&\sup_{x_I\in\bR^{|I|}}\sup_{y_I\in\bR^{|I|}}\left|L_{(h_I)}(x_I-y_{I})\right|\leq\left\|L_{(h_I)}\right\|_{\infty}\leq C_I(\mathbf{K},q)\prod_{j\in I}h_j^{-\lambda_j-1}<\infty,
\\*[2mm]
&&\ C_I(\mathbf{K},q):=\frac{\mathbf{A}}{(2\pi)^{\frac{|I|}{2}}}\left\{\left\|\widehat{K_I}g_I\right\|_{2}\vee\left\|\widehat{K_I}g_I\right\|_{1}\vee\left(\max_{j\in I}\left\|D_j^1\widehat{K_I}g_I\right\|_{1}\right)\vee\left\|\widehat{K_I}\varphi_I\right\|_{2}\vee\left\|\widehat{K_I}\varphi_I\right\|_{1}\right\},\nonumber
\end{eqnarray}
where $\varphi_I(t_I):=\sup_{j\in I}\left|t_j\right|g_I(t_I)$ and $g_I$ is given in (\ref{eq:auxiliaryfunction});
\begin{eqnarray}
\label{eq:variancebound1}
\sup_{x_I\in\bR^{|I|}}\left(\bE\left|L_{(h_I)}(x_I-Y_{1,I})\right|^2\right)^{\frac{1}{2}}\leq\sqrt{\left\|f_I\star q_I\right\|_{\infty}}\left\|L_{(h_I)}\right\|_2 \leq \sqrt{\left\|q_I\right\|_{\infty}}C_I(\mathbf{K},q)\prod_{j\in I}h_j^{-\lambda_j-\frac{1}{2}}.
\end{eqnarray}

Next, set $x_I$ and $\overline{x}_I$ be arbitrary fixed in $\bR^{|I|}$. For any $t_I\in\bR^{|I|}$
\begin{eqnarray*}
\left|e^{-i<t_I,x_I>}-e^{-i<t_I,\overline{x}_I>}\right|&=&\left|\prod_{j\in I}e^{-it_jx_j}-\prod_{j\in I}e^{-it_j\overline{x}_j}\right|
\\*[2mm]
&\leq&|I|\sup_{j\in I}\left|e^{-it_jx_j}-e^{-it_j\overline{x}_j}\right|
\\*[2mm]
&\leq&|I|\sup_{j\in I}\left|t_j\right|\sup_{j\in I}\left|x_j-\overline{x}_j\right|.
\end{eqnarray*}
Therefore, for any $y_I\in\bR^{|I|}$
\begin{eqnarray}
\label{eq:uniformbound12}
&&\left|L_{(h_I)}(x_I-y_{I})-L_{(h_I)}(\overline{x}_I-y_{I})\right|\nonumber
\\*[2mm]
&&\leq\frac{1}{(2\pi)^{|I|}}\int_{\bR^{|I|}}\left|\frac{\widehat{K_{h_I}}(t_I)}{\widehat{q_I}(t_I)}\right|\left|e^{-i<t_I,x_I>}-e^{-i<t_I,\overline{x}_I>}\right|\rd t_I\nonumber
\\*[2mm]
&&\leq n|I|C_I(\mathbf{K},q)\prod_{j\in I}h_j^{-\lambda_j-1}\sup_{j\in I}\left|x_j-\overline{x}_j\right|;
\end{eqnarray}
\begin{eqnarray}
\label{eq:variancebound2}
&&\left(\bE\left|L_{(h_I)}(x_I-Y_{1,I})-L_{(h_I)}(\overline{x}_I-Y_{1,I})\right|^2\right)^{\frac{1}{2}}\nonumber
\\*[2mm]
&&\leq\left(\frac{\left\|f_I\star q_I\right\|_{\infty}}{(2\pi)^{|I|}}\int_{\bR^{|I|}}\left|\frac{\widehat{K_{h_I}}(t_I)}{\widehat{q_I}(t_I)}\right|^2\left|e^{-i<t_I,x_I>}-e^{-i<t_I,\overline{x}_I>}\right|^2\rd t_I\right)^{\frac{1}{2}}\nonumber
\\*[2mm]
&&\leq n|I|\sqrt{1\vee\left\|q_I\right\|_{\infty}}C_I(\mathbf{K},q)\prod_{j\in I}h_j^{-\lambda_j-\frac{1}{2}}\sup_{j\in I}\left|x_j-\overline{x}_j\right|;
\end{eqnarray}

Consider now the normalized empirical process
$$
\overline{\xi}_{h_I}(x_I):=\left(C_I(\mathbf{K},q)\sqrt{\frac{2(1\vee\left\|q_I\right\|_{\infty})}{n\prod_{j\in I}h_j^{2\lambda_j+1}}}\right)^{-1}\xi_{h_I}(x_I).
$$

In view of Bernstein inequality, (\ref{eq:uniformbound12}), (\ref{eq:variancebound1}), (\ref{eq:uniformbound12}) and (\ref{eq:variancebound2}), $\forall z>0$,
\begin{eqnarray}
\label{eq:bernstein1}
\bP\left\{\left|\overline{\xi}_{h_I}(x_I)\right|>z\right\}&\leq& 2\exp\left\{-\frac{z^2}{A^2(x_I)+zB(x_I)}\right\};
\\*[2mm]
\label{eq:bernstein2}
\bP\left\{\left|\overline{\xi}_{h_I}(x_I)-\overline{\xi}_{h_I}(\overline{x}_I)\right|>z\right\}&\leq& 2\exp\left\{-\frac{z^2}{\ra^2(x_I,\overline{x}_I)+z\rb(x_I,\overline{x}_I)}\right\},
\end{eqnarray}
where $A(x_I):=1$, $B(x_I):=\left(n\prod_{j\in I}h_j^{2\lambda_j+1}\right)^{-\frac{1}{2}}\leq 1$ and
\begin{equation}
\label{eq:semimetrics}
\ra(x_I,\overline{x}_I)=\rb(x_I,\overline{x}_I):=2\wedge\left\{n|I|\sup_{j\in I}\left|x_j-\overline{x}_j\right|\right\}.
\end{equation}
It is easily seen that $\ra(\cdot,\cdot)$ is a semi-metric on $\bR^{|I|}$.

\paragraph{2) Supremum-norm over totally bounded sets: } In this step we obtain bounds of the supremum-norm of the normalized empirical process $\overline{\xi}_{h_I}(\cdot)$ over totally bounded sets by applying Proposition 1 in Lepski \cite{lepski:upperfunctions} with $\mT=\bR^{|I|}$, $\mS=\bR$, $\chi=\overline{\xi}_{h_I}$ and $\Psi(\cdot)=\left|\;\cdot\;\right|$. Then we have to check Assumptions 1, 2 and 3 required in the latter 
Proposition and to match the notations used in the present paper and in Lepski \cite{lepski:upperfunctions}.

Note first that, in view of (\ref{eq:bernstein1}), (\ref{eq:bernstein2}) and (\ref{eq:semimetrics}), Assumption 1 is fulfilled with $c=2$. Next, consider the family of closed balls
$$
\bB_{\frac{R}{2}}(t_I):=\left\{x_I\in\bR^{|I|}:\;\sup_{j\in I}\left|x_j-t_j\right|\leq R/2\;\right\},\quad R\geq 1,\; t_I\in\bR^{|I|}.
$$
In view of the continuity property of the Fourier transforms and the definition of the semi-metrics $\ra$ and $\rb$, it is obvious that Assumption 2 is also satisfied with $\Theta=\bB_{\frac{R}{2}}(t_I)$ , $\overline{A}_\Theta=1$ and $\overline{B}_\Theta=\left(n\prod_{j\in I}h_j^{2\lambda_j+1}\right)^{-\frac{1}{2}}$.

\smallskip

Let $s:\;\bR\rightarrow\bR_+\backslash\{0\}$ defined by $s(z):=(0,01+z^8)^{-1}$. Obviously $\sum_{k\geq0}s\big(2^{k/2}\big)\leq 1$ and, for any $z>0$,
\begin{equation}
\label{eq:entropy1}
\mE_{\Theta,\ra}\left(z(48\delta)^{-1}s(\delta)\right)\leq|I|\left[\ln\left(\frac{Rn|I|}{z(48\delta)^{-1}s(\delta)}\right)\right]_+,\quad\forall\delta>0,
\end{equation}
where $\mE_{\Theta,\ra}\left(\delta\right)$, $\delta>0$, denotes the entropy of $\Theta$ measured in $\ra$. Then, for any $z>0$, there exists $\delta_*>0$ small enough such that
\begin{eqnarray*} e_s^{(\ra)}(z,\Theta)&:=&\sup_{\delta>0}\delta^{-2}\mE_{\Theta,\ra}\left(z(48\delta)^{-1}s(\delta)\right)=\sup_{\delta>\delta_*}\delta^{-2}\mE_{\Theta,\ra}\left(z(48\delta)^{-1}s(\delta)\right)<\infty;
\\*[2mm] e_s^{(\rb)}(z,\Theta)&:=&\sup_{\delta>0}\delta^{-1}\mE_{\Theta,\rb}\left(z(48\delta)^{-1}s(\delta)\right)=\sup_{\delta>\delta_*}\delta^{-1}\mE_{\Theta,\rb}\left(z(48\delta)^{-1}s(\delta)\right)<\infty.
\end{eqnarray*}
Thus, Assumption 3 in Lepski \cite{lepski:upperfunctions} is fulfilled and Proposition 1 in the latter paper can be applied. Let us compute the quantities which appear in this result.

\smallskip

Choose $\vec{s}=(s,s)$, $\varkappa=(2\overline{A}_\Theta,2\overline{B}_\Theta)$ and $\e=\sqrt{2}-1$. Since $\overline{A}_\Theta\vee\overline{B}_\Theta\leq 1$ and $\ra(x_I,\overline{x}_I)=\rb(x_I,\overline{x}_I)\leq 2$, $\forall x_I,\overline{x}_I\in\bR^{|I|}$, we straightforwardly get
\begin{eqnarray*}
e_{\vec{s}}(\varkappa,\Theta)&:=&e_s^{(\ra)}(2\overline{A}_\Theta,\Theta)+e_s^{(\rb)}(2\overline{B}_\Theta,\Theta)
\\*[2mm] &\leq&\sup_{\delta>0,61}\delta^{-2}\mE_{\Theta,\ra}\left(2(48\delta)^{-1}s(\delta)\right)+\sup_{\delta>0,61}\delta^{-1}\mE_{\Theta,\rb}\left(2(48\delta)^{-1}s(\delta)\right)
\\*[2mm] 
&\leq& 4,5|I|\left[\ln\left(Rn|I|\right)\right]_++8,5;
\end{eqnarray*}
\begin{eqnarray*}
U_{\vec{s}}^{(\e)}(y,\varkappa,\Theta)&:=&\varkappa_1\sqrt{2[1+\e^{-1}]^2e_{\vec{s}}(\varkappa,\Theta)+y}+\varkappa_2\left(2[1+\e^{-1}]^2e_{\vec{s}}(\varkappa,\Theta)+y\right)
\\*[2mm] 
&\leq& 2\sqrt{31|I|\ln\left(Rn|I|\right)+59+y}+\frac{2\left(31|I|\ln\left(Rn|I|\right)+59+y\right)}{\sqrt{n\prod_{j\in I}h_j^{2\lambda_j+1}}}.
\end{eqnarray*}

Thus, it follows from Proposition 1 in Lepski \cite{lepski:upperfunctions} that, for any $y\geq 1$ and any $\mathbf{r}\geq 1$,
\begin{equation}
\label{eq:upperbound1}
\bE\left\{\sup_{x_I\in\bB_{\frac{R}{2}}(t_I)}\left|\overline{\xi}_{h_I}(x_I)\right|-U_{\vec{s}}^{(\e)}(y,\varkappa,\Theta)\right\}_+^r\leq 4\Gamma(\mathbf{r}+1)\left[2y^{-1}U_{\vec{s}}^{(\e)}(y,\varkappa,\Theta)\right]^{\mathbf{r}}e^{-\frac{y}{2}}.
\end{equation}

\paragraph{3) Supremum-norm over the whole space: } Let $x_I\in\bR^{|I|}$ be arbitrary fixed and $y_I\in\bR^{|I|}$ be such that $\sup_{j\in I}|x_j-y_j|\geq n$. By integration by parts, we easily get
\begin{equation}
\label{eq:separation}
\left|L_{(h_I)}(x_I-y_I)\right|\leq\frac{\max_{j\in I}\left\|D_j^1\left(\widehat{K_{h_I}}/\widehat{q_I}\right)\right\|_1}{(2\pi)^{|I|}\sup_{j\in I}\left|x_j-y_j\right|}\leq \frac{C_I(\mathbf{K},q)}{n\prod_{j\in I}h_j^{\lambda_j+1}}\leq\frac{C_I(\mathbf{K},q)}{n\prod_{j\in I}h_j^{2\lambda_j+1}},
\end{equation}
in view of Assumption (N2) on the errors.

Consider the collection of closed balls $\Big\{\bB_{\frac{n}{2}}(n\textbf{j}),\;\textbf{j}\in\bZ^{|I|}\Big\}$. Obviously this collection is a countable cover of $\bR^{|I|}$. Put, for any $\textbf{j}\in\bZ^{|I|}$,
$$
\rf_{\textbf{j}}:=\int_{\bB(\textbf{j})}f_I\star q_I(x_I)\rd x_I,\qquad\bB(\textbf{j}):=\bigcup_{\textbf{k}\in\bZ^{|I|}:\;\bB_{\frac{n}{2}}(n\textbf{j})\cap\bB_{\frac{n}{2}}(n\textbf{k})\neq\emptyset}\bB_{\frac{n}{2}}(n\textbf{k}).
$$
It is easily checked that
\begin{equation}
\label{eq:probabilitybound1}
\sum_{\textbf{j}\in\bZ^{|I|}}\rf_{\textbf{j}}=\int_{\bR^{|I|}}f_I\star q_I(x_I)\left[\sum_{\textbf{j}\in\bZ^{|I|}}\textbf{1}_{\bB(\textbf{j})}(x_I)\right]\rd x_I\leq 4^{|I|}.
\end{equation}

\smallskip

Set $\textbf{j}\in\bZ^{|I|}$ such that $\rf_{\textbf{j}}\geq n^{-v}$, where $v\geq 1$ is specified later. If $y=2\ln(1/\rf_{\textbf{j}})+(\mathbf{r}+1)\ln(n)$, we get from (\ref{eq:upperbound1})
\begin{equation*}
\bE\left\{\sup_{x_I\in\bB_{\frac{n}{2}}(n\textbf{j})}\left|\xi_{h_I}(x_I)\right|-\gamma_{\infty,I}^{(v)}(\mathbf{r})\sqrt{\frac{\ln(n)}{n\prod_{j\in I}h_j^{2\lambda_j+1}}}\right\}_+^{\mathbf{r}}\leq 2^{\mathbf{r}+2}\Gamma(\mathbf{r}+1)\left[\gamma_{\infty,I}^{(v)}(\mathbf{r})\right]^{\mathbf{r}}\rf_{\textbf{j}}n^{-\frac{\mathbf{r}+1}{2}},
\end{equation*}
where $\gamma_{\infty,I}^{(v)}(\mathbf{r}):=4C_I(\mathbf{K},q)\sqrt{2(1\vee\left\|q_I\right\|_{\infty})}(93|I|\ln(|I|)+60+2v+\mathbf{r})$, since $n\prod_{j\in I}h_j^{2\lambda_j+1}\geq\ln(n)$.

\smallskip

Thus, in view of (\ref{eq:probabilitybound1}), we obtain
\begin{equation}
\label{eq:upperbound2}
\bE\left\{\sup_{x_I\in\Theta_1}\left|\xi_{h_I}(x_I)\right|-\gamma_{\infty,I}^{(v)}(\mathbf{r})\sqrt{\frac{\ln(n)}{n\prod_{j\in I}h_j^{2\lambda_j+1}}}\right\}_+^{\mathbf{r}}\leq 2^{\mathbf{r}+2+2|I|}\Gamma(\mathbf{r}+1)\left[\gamma_{\infty,I}^{(v)}(\mathbf{r})\right]^{\mathbf{r}} n^{-\frac{\mathbf{r}+1}{2}},
\end{equation}
where $\Theta_1:=\cup_{\textbf{j}\in\bZ^{|I|}:\rf_{\textbf{j}}\geq n^{-v}}\bB_{\frac{n}{2}}(n\textbf{j})$.

\smallskip

Set $\textbf{j}\in\bZ^{|I|}$ such that $\rf_{\textbf{j}}< n^{-v}$ and $x_I\in \bB_{\frac{n}{2}}(n\textbf{j})$. In view of (\ref{eq:uniformbound1}) and (\ref{eq:separation}) we get, for any $k=1,\ldots,n$,
\begin{eqnarray}
\label{eq:expectation}
\bE\left|L_{(h_I)}(x_I-Y_{k,I})\right|&=&\bE\left\{\left|L_{(h_I)}(x_I-Y_{k,I})\right|\textbf{1}_{\bB(\textbf{j})}(Y_{k,I})\right\}\nonumber
\\*[2mm]
&&+\bE\left\{\left|L_{(h_I)}(x_I-Y_{k,I})\right|\textbf{1}_{\bR^{|I|}\backslash\bB(\textbf{j})}(Y_{k,I})\right\}\nonumber
\\*[2mm]
&\leq&\bP\left\{Y_{k,I}\in\bB(\textbf{j})\right\}\frac{C_I(\mathbf{K},q)}{\prod_{j\in I}h_j^{2\lambda_j+1}}+\frac{C_I(\mathbf{K},q)}{n\prod_{j\in I}h_j^{2\lambda_j+1}}\nonumber
\\*[2mm]
&\leq&\frac{2C_I(\mathbf{K},q)}{n\prod_{j\in I}h_j^{2\lambda_j+1}},
\end{eqnarray}
since $\rf_{\textbf{j}}:=\bP\left\{Y_{k,I}\in\bB(\textbf{j})\right\}\leq n^{-v}$, $v\geq 1$ and $\sup_{j\in I}|x_j-Y_{k,j}|\geq n$ when $Y_{k,I}\in \bR^{|I|}\backslash\bB(\textbf{j})$.

\smallskip

Introduce random events
$$
D_{\textbf{j}}:=\left\{\sum_{k=1}^n\textbf{1}_{\bB(\textbf{j})}(Y_{k,I})\geq 2\right\},\;\textbf{j}\in\bZ^{|I|},\quad D:=\bigcup_{\textbf{j}\in\bZ^{|I|}:\rf_{\textbf{j}}< n^{-v}}D_{\textbf{j}}.
$$
Let $\overline{D}$ be the complementary to $D$. If $\overline{D}$ holds then, in view of (\ref{eq:uniformbound1}) and (\ref{eq:separation}),
\begin{equation}
\label{eq:empiricalbound}
n^{-1}\sum_{k=1}^n\left|L_{(h_I)}(x_I-Y_{k,I})\right|\leq\frac{2C_I(\mathbf{K},q)}{n\prod_{j\in I}h_j^{2\lambda_j+1}},\quad\forall x_I\in\Theta_2:=\bR^{|I|}\backslash\Theta_1.
\end{equation}
Since $n\prod_{j\in I}h_j^{2\lambda_j+1}\geq\ln(n)$, we get from (\ref{eq:expectation}) and (\ref{eq:empiricalbound})
$$
\sup_{x_I\in\Theta_2}\left|\xi_{h_I}(x_I)\right|\textbf{1}_{\overline{D}}\leq\gamma_{\infty,I}^{(v)}(\mathbf{r})\sqrt{\frac{\ln(n)}{n\prod_{j\in I}h_j^{2\lambda_j+1}}}
$$
and, taking into account that $\sup_{x_I\in\Theta_2}\left|\xi_{h_I}(x_I)\right|\leq 2C_I(\mathbf{K,q})n$,
\begin{equation}
\label{eq:upperbound3}
\bE\left\{\sup_{x_I\in\Theta_2}\left|\xi_{h_I}(x_I)\right|-\gamma_{\infty,I}^{(v)}(\mathbf{r})\sqrt{\frac{\ln(n)}{n\prod_{j\in I}h_j^{2\lambda_j+1}}}\right\}_+^{\mathbf{r}}\leq \left[2C_I(\mathbf{K,q})\right]^{\mathbf{r}} n^{\mathbf{r}}\bP(D).
\end{equation}

\smallskip

Let $\textbf{j}\in\bZ^{|I|}$ satisfying $\rf_{\textbf{j}}< n^{-v}$ be arbitrary fixed. In view of Markov inequality one has for any $z>0$
$$
\bP(D_\textbf{j})\leq e^{-2z}\left[\bE\left\{e^{z\textbf{1}_{\bB(\textbf{j})}(Y_{1,I})}\right\}\right]^n\leq\exp\left\{-2z+n(e^z-1)\rf_{\textbf{j}}\right\},
$$
since the $Y_{k,I}$'s are i.i.d. random vectors. Minimizing the right hand side in $z>0$ we obtain
\begin{equation}
\label{eq:probabilitybound2}
\bP(D_\textbf{j})\leq(e/2)^2(n\rf_{\textbf{j}})^2\leq 2\rf_{\textbf{j}}n^{2-v}.
\end{equation}

\smallskip

Thus, choosing $v=1,5\mathbf{r}+2,5$, it follows from (\ref{eq:probabilitybound1}), (\ref{eq:upperbound2}), (\ref{eq:upperbound3}) and (\ref{eq:probabilitybound2})
\begin{equation}
\label{eq:upperbound4}
\bE\left\{\left\|\xi_{h_I}\right\|_{\infty}-\gamma_{\infty,I}(\mathbf{r})\sqrt{\frac{\ln(n)}{n\prod_{j\in I}h_j^{2\lambda_j+1}}}\right\}_+^{\mathbf{r}}\leq 2^{\mathbf{r}+3+2|I|}\Gamma(\mathbf{r}+1)\left[\gamma_{\infty,I}(\mathbf{r})\right]^{\mathbf{r}} n^{-\frac{\mathbf{r}+1}{2}},
\end{equation}
where $\gamma_{\infty,I}(\mathbf{r}):=\gamma_{\infty,I}^{(1,5\mathbf{r}+2,5)}(\mathbf{r})$.

\smallskip

Finally, in view of the definition of $\cH_{\infty,I}$,
\begin{eqnarray}
\label{eq:upperbound5}
&&\left\{\bE\sup_{h_I\in\cH_{\infty,I}}\left[\left\|\xi_{h_I}\right\|_{\infty}-\gamma_{\infty,I}(\mathbf{r})\cU_{\infty}(h_I)\right]_+^{\mathbf{r}}\right\}^{\frac{1}{\mathbf{r}}}\leq c_{\infty}(\mathbf{r})n^{-\frac{1}{2}},
\\*[2mm]
&& c_{\infty}(\mathbf{r}):=\left[\Gamma(\mathbf{r}+1)\right]^{\frac{1}{\mathbf{r}}}\sup_{n\in\bN^*}\sup_{I\in\cI_d^{\diamond}}\left\{\gamma_{\infty,I}(\mathbf{r})\left[2^{\mathbf{r}+3+2|I|}\right]^{\frac{1}{\mathbf{r}}}[\log_2(n)]^{\frac{|I|}{\mathbf{r}}}n^{-\frac{1}{2\mathbf{r}}}\right\}<\infty.\nonumber
\end{eqnarray}
\epr

\subsection{Proof of Lemma \ref{lem:empiricalupperbound1}} 

Assume that $\mP\neq\big\{\overline{\emptyset}\big\}$. Set $f\in\bF_p\left[\;\mP\;\right]$ and let $\mathbf{r}\in\{\mathbf{r}_1,\mathbf{r}_2,\mathbf{r}_4\}$ be arbitrary fixed. We obtain Lemma \ref{lem:empiricalupperbound1} by applying Proposition \ref{prop:empiricalupperbound1}. We divide this proof into two steps.

\paragraph{1)} Note that 
\begin{equation*}
\label{eq:lemma11}
\xi_{p}\leq\sum_{I\in\cI_d^{\diamond}}\sup_{h_I\in\cH_{p,I}}\left[\left\|\xi_{h_I}\right\|_{p}-\gamma_{p,I}(\mathbf{r})\cU_p(h_I)\;\right]_+,
\end{equation*}
since $\gamma_{p,I}(\mathbf{r})$ increase with $\mathbf{r}$. In view of Proposition \ref{prop:empiricalupperbound1}, if $p\in(1,+\infty]$ and $n\geq 3$,
\begin{eqnarray*}
\label{eq:lemma12}
\left(\bE_f\left|\xi_{p}\right|^{\mathbf{r}}\right)^{\frac{1}{\mathbf{r}}}\leq\textbf{c}_{p,1}(\mathbf{r})n^{-\frac{1}{2}},\quad\textbf{c}_{p,1}(\mathbf{r}):=d|\mP|^2c_p(\mathbf{r}).
\end{eqnarray*}

\paragraph{2)} For any $p\geq 1$
\begin{eqnarray*}
\overline{G}_{p}&\leq&  1+\left\|\textbf{K}\right\|_1^d\sup_{I\in\cI_d^{\diamond}}\sup_{h_I\in\cH_{p,I}}\left\{\left[\left\|\xi_{h_I}\right\|_{p}-\overline{\gamma}_p\cU_p(h_I)\;\right]_++\overline{\gamma}_p\cU_p(h_I)+\left\|\bE_f\left\{\widetilde{f}_{h_I}\right\}\right\|_{p}\right\}
\\*[2mm]
&\leq& 1+\left\|\textbf{K}\right\|_1^{d}\left(\overline{\xi}_p+\overline{\gamma}_p\overline{\cU}_p+\left\|\textbf{K}\right\|_1^{d}\mathbf{f}_p\right),
\\*[2mm]
\overline{\gamma}_p&:=&\sup_{I\in\cI_d^{\diamond}}\gamma_{p,I}(\mathbf{r}_4\md^2),\quad\overline{\xi}_p:=\sup_{I\in\cI_d^{\diamond}}\sup_{h_I\in\cH_{p,I}}\left[\left\|\xi_{h_I}\right\|_{p}-\gamma_{p,I}(\mathbf{r}_4\md^2)\cU_p(h_I)\;\right]_+;
\\*[2mm]
\overline{\mathbf{f}}_{p}&\leq& \md^2\left\|\textbf{K}\right\|_1^{d}\Bigg[\overline{G}_{p}+\left\|\textbf{K}\right\|_1^{d}\mathbf{f}_p\Bigg]^{\md^2-1}\leq\md^2\left\|\textbf{K}\right\|_1^{d}\Bigg[1+\left\|\textbf{K}\right\|_1^{d}\left(\overline{\xi}_p+\overline{\gamma}_p\overline{\cU}_p+\left\|\textbf{K}\right\|_1^{d}\mathbf{f}_p+\mathbf{f}_p\right)\Bigg]^{\md^2}
\end{eqnarray*}

\smallskip

Below we use the following trivial equality: for any random variable $Y$
\begin{equation}
\label{eq:lemma17}
\left(\bE_f\left|Y^{\md^2}\right|^{\mathbf{r}}\right)^{\frac{1}{\mathbf{r}}}=\Bigg[\left(\bE_f\left|Y\right|^{\mathbf{r}\md^2}\right)^{\frac{1}{\mathbf{r}\md^2}}
\Bigg]^{\md^2},
\end{equation}

In view of Proposition \ref{prop:empiricalupperbound1}, if $p\in(1,+\infty]$ and $n\geq 3$, $\left(\bE_f\left|\overline{\mathbf{f}}_{p}\right|^{\mathbf{r}}\right)^{\frac{1}{\mathbf{r}}}\leq \textbf{c}_{p,2}(\mathbf{r},\mathbf{f}_p)$ with 
$$
\textbf{c}_{p,2}(\mathbf{r},\mathbf{f}_p):=\md^2\left\|\textbf{K}\right\|_1^{d}\left[1+\left\|\textbf{K}\right\|_1^{d}\left(d|\mP|^2c_p(\mathbf{r}\md^2)+\overline{\gamma}_p\overline{\cU}_p+\left\|\textbf{K}\right\|_1^{d}\mathbf{f}_p+\mathbf{f}_p\right)\right]^{\md^2}.
$$
Thus, we finish the proof of Lemma \ref{lem:empiricalupperbound1}.
\epr

\bibliographystyle{agsm}

\end{document}